\documentclass[letterpaper, 12 pt]{article}
\usepackage{amsthm, amsmath, amssymb, amsfonts, graphicx}

\newcommand{\setleftmargin}[1]{
        \addtolength{\textwidth}{\oddsidemargin}
        \addtolength{\textwidth}{1in}
        \addtolength{\textwidth}{-#1}
        \setlength{\oddsidemargin}{-1in}
        \addtolength{\oddsidemargin}{#1}
        \setlength{\evensidemargin}{\oddsidemargin}
}
\newcommand{\setrightmargin}[1]{
        \setlength{\textwidth}{8.5in}
        \addtolength{\textwidth}{-\oddsidemargin}
        \addtolength{\textwidth}{-1in}
        \addtolength{\textwidth}{-#1}
}
\newcommand{\settopmargin}[1]{ 
        \addtolength{\textheight}{\topmargin}
        \addtolength{\textheight}{1in}
        \addtolength{\textheight}{\headheight}
        \addtolength{\textheight}{\headsep}
        \addtolength{\textheight}{-#1}
        \setlength{\topmargin}{-1in}
        \addtolength{\topmargin}{-\headheight}
        \addtolength{\topmargin}{-\headsep}
        \addtolength{\topmargin}{#1}
}
\newcommand{\setbottommargin}[1]{
        \setlength{\textheight}{11in}
        \addtolength{\textheight}{-\topmargin}
        \addtolength{\textheight}{-1in}
        \addtolength{\textheight}{-\footskip}
        \addtolength{\textheight}{-#1}
}
\newcommand{\setallmargins}[1]{
        \settopmargin{#1}
        \setbottommargin{#1}
        \setleftmargin{#1}
        \setrightmargin{#1}
}

\setallmargins{1 in}

\newcommand{\Flag}{\mathrm{Flag}}

\newcommand{\ZZ}{\mathbb{Z}}
\newcommand{\QQ}{\mathbb{Q}}
\newcommand{\RR}{\mathbb{R}}
\newcommand{\CC}{\mathbb{C}}
\newcommand{\PP}{\mathbb{P}}

\newcommand{\FF}{\mathcal{F}}

\newcommand{\GG}{\mathcal{G}}

\newcommand{\gen}{\mathrm{gen}}
\newcommand{\Poly}{\mathrm{Poly}}

\newcommand{\into}{\hookrightarrow}
\newcommand{\onto}{\twoheadrightarrow}

\newcommand{\ra}{\rightarrow}

\DeclareMathOperator{\Span}{Span}

\DeclareMathOperator{\diag}{diag}

\DeclareMathOperator{\Unif}{Uniform}

\DeclareMathOperator{\Spec}{Spec}
\DeclareMathOperator{\Proj}{Proj}

\newcommand{\Ko}{\mathrm{K}^{\circ}}

\newcommand{\PGL}{\mathrm{PGL}}

\newcommand{\OO}{\mathcal{O}}

\newcommand{\DD}{\mathcal{D}}
\newcommand{\EE}{\mathcal{E}}
\newcommand{\oDD}{\mathring{\DD}}

\newcommand{\sym}{\sim}

\newcommand{\Go}{\mathring{G}}
\newcommand{\Xo}{\mathring{X}}
\newcommand{\TorTor}{\mathcal{T}or}

\newcommand{\isomorph}{\cong}
\newtheorem{conj}{Conjecture}[section]

\newtheorem{theorem}[conj]{Theorem}
\newtheorem{prop}[conj]{Proposition}

\newtheorem{lemma}[conj]{Lemma}

\newtheorem{cor}[conj]{Corollary}

\newtheorem*{f-vector*}{The $f$-Vector Conjecture}
\newtheorem*{f-vector-sp*}{The $f$-Vector Conjecture, Continued}
\newtheorem*{MT}{Main Theorem}
\newtheorem*{PF}{Proposition \ref{ProductFormula}}
\newtheorem*{PT}{Proposition \ref{Positivity}}

\theoremstyle{definition}

\begin{document}

\title{A Matroid Invariant \emph{via} the $K$-Theory of the Grassmannian}
\author{David E Speyer \\ \texttt{speyer@umich.edu}}
\date{}
\maketitle

\abstract{Let $G(d,n)$ denote the Grassmannian of $d$-planes in $\CC^n$ and let $T$ be the torus $(\CC^*)^n/\diag(\CC^*)$ which acts on $G(d,n)$. Let $x$ be a point of $G(d,n)$ and let $\overline{Tx}$ be the closure of the $T$-orbit through $x$. Then the class of the structure sheaf of $\overline{Tx}$ in the $K$-theory of $G(d,n)$ depends only on which Pl\"ucker coordinates of $x$ are nonzero -- combinatorial data known as the matroid of $x$. In this paper, we will define a certain map of additive groups from $K^{\circ}(G(d,n))$ to $\ZZ[t]$. Letting $g_x(t)$ denote the image of $(-1)^{n-\dim Tx} [ \OO_{\overline{Tx}}]$, $g_x$ behaves nicely under the standard constructions of matroid theory. Specifically, $g_{x_1 \oplus x_2}(t)=g_{x_1}(t) g_{x_2}(t)$, $g_{x_1 +_2 \ x_2}(t)=g_{x_1}(t) g_{x_2}(t)/t$, $g_x(t) = g_{x^{\perp}}(t)$ and $g_x$ is unaltered by series and parallel extensions. Furthermore, the coefficients of $g_x$ are nonnegative. The existence of this map implies bounds on (essentially equivalently) the complexity of Kapranov's Lie complexes \cite{CQoG1}, Hacking, Keel and Tevelev's very stable pairs \cite{HKT} and the author's tropical linear spaces when they are realizable in characteristic zero \cite{TLS}. Namely, in characteristic zero, a Lie complex or the underlying $d-1$ dimensional scheme of a very stable pair can have at most $\frac{(n-i-1)!}{(d-i)!(n-d-i)!(i-1)!}$ strata of dimensions $n-i$ and $d-i$ respectively and a tropical linear space realizable in characteristic zero can have at most this many $i$-dimensional bounded faces.}

\section{Motivation and Introduction} \label{Intro}

Let $K= \bigcup_{n=1}^{\infty} \CC((t^{1/n}))$, the field of Puiseux series, and let $v: K^* \to \QQ$ be the map which assigns to a power series its order of vanishing; in other words, if $x=\sum_{i \geq M} a_i t^{i/N}$ with $a_M \neq 0$ then $v(x)=M/N$. Suppose that we have a $K$-valued point of the Grassmannian $G(d,n)$ none of whose Pl\"ucker coordinates $p_I(x)$ are zero. In \cite{TTG}, Sturmfels and I attempted to determine the possible $\binom{n}{d}$-tuples of rational numbers $v(p_I(x))$.  

The $p_I(x)$ obey the Pl\"ucker relations:
$$p_{Sij} p_{Sk \ell} - p_{Sik} p_{Sj \ell} + p_{Si \ell} p_{Sjk}=0$$
for any $S \in \binom{[n]}{d-2}$ \footnote{Through out this paper, we use the following combinatorial notations: $[n] = \{ 1,2,\ldots, n \}$ and, for any set $S$, $\binom{S}{d}$ is the set of $d$ element subsets of $S$. We will sometimes, as in the equation above, use the shorthand $Si$ for $S \cup \{ i \}$.} and $i<j<k<\ell$ in $[n] \setminus S$. As a consequence, we deduce that among the three numbers 
$$v(p_{Sij})+v(p_{Sk \ell}) \quad v(p_{Sik})+v(p_{Sj \ell}) \quad v(p_{Si \ell})+v(p_{Sjk})$$
the minimum occurs at least twice. An $\binom{n}{d}$-tuple of rational numbers $P_I$ obeying this condition is a valuated matroid in the sense of Dress and Wenzel \cite{DW} or, in the terminology I suggested in \cite{TLS}, a tropical Pl\"ucker vector. In \cite{TLS}, I explained how to associate to a tropical Pl\"ucker vector a $d$-dimensional polyhedral complex in $\RR^n$, called a tropical linear space, and showed how to perform operations on tropical linear spaces such as orthogonal complement, intersection and span in a manner analogous to those operations for ordinary linear spaces. In that paper, it was shown that any tropical linear space built out of repeatedly applying these simple operations had the same number of $c$-dimensional bounded faces -- specifically, $\frac{(n-c-1)!}{(d-c)!(n-d-c)!(c-1)!}$. (This is a slight simplification -- see \cite{TLS} for details.) It was conjectured in that paper that this was the maximal number of $c$-dimensional bounded faces for any tropical linear space for given $(d,n)$. My motivation in beginning the research reported in this paper was to prove this conjecture; in this paper we will give a proof for tropical Pl\"ucker vectors which arise as $v(p_I(K))$ as above.

There is a simple polyhedral construction that lets us understand the combinatorial meaning of being a tropical Pl\"ucker vector. Let $\Delta(d,n)$ denote the $(d,n)$-hypersimplex, the convex hull of the $\binom{n}{d}$ points $e_{i_1}+\ldots+e_{i_d}$ where $(i_1,\ldots,i_d)$ runs over $\binom{[n]}{d}$. Let $I \mapsto P_I$ be a function $\binom{[n]}{d} \to \QQ$. We define a polyhedral subdivision of $\Delta(d,n)$ as follows: Let $Q$ denote the convex hull in $\Delta(d,n) \times \RR$ of the points $(e_{i_1}+\ldots + e_{i_d}, P_{i_1 \ldots i_d})$ and let $Q + \RR_{\geq 0}$ denote the Minkowski sum of $Q$ with $\{ 0 \} \times \RR_{\geq 0} \subset \RR^n \times \RR$. Take the facets of $Q + \RR_{\geq 0}$ whose outward pointing normal vectors have negative components in the last coordinate and project them down to $\Delta(d,n)$. This gives a polyhedral subdivision $\DD_P$ of $\Delta(d,n)$ which is known as the regular subdivision associated to $P$. Then $P$ is a tropical Pl\"ucker vector if and only if, for every face $F$ of $\Delta$, the vertices of $F$ form the bases of a matroid (when considered as a subset of $\binom{[n]}{d}$.) \footnote{For all matroid related terminology in this paper, our references are \cite{White1} and \cite{White2}. We will use very little deep material from the theory of matroids -- worried readers are urged to see the note at the end of this section. } When this condition holds,  the number of $c$-dimensional bounded faces of the corresponding tropical linear space is the number of $(n-c)$-dimensional interior faces in $\DD$. See \cite{TLS} for more on this construction and see \cite{Zieg}, chapter 5, for more material on regular subdivisions in general. We can now state precisely the $f$-vector conjecture from \cite{TLS}.

\begin{f-vector*}
Let $P$ be a tropical Pl\"ucker vector and define $\DD_P$ as above. Then $\DD_P$ has at most $\frac{(n-c-1)!}{(d-c)! (n-d-c)! (c-1)!}$ interior faces of dimension $n-c$, with equality if and only if all of the facets of $\DD_p$ correspond to series-parallel matroids.
\end{f-vector*}

Series-parallel matroids are a certain well known class of matroid, we will give the definition of a series-parallel matroid in section \ref{GammaConstruction}.

In the case where the $P_I$ are $v(p_I(x))$ for some $x \in G(d,n)(K)$ there is a geometric meaning to the decomposition $\DD_{P}$. The $(n-1)$-dimensional torus $T=(\CC^*)^n/\diag(\CC^*)$ acts on $G(d,n)$ and $\overline{T x}$ is the toric variety over $K$ associated to $\Delta(d,n)$. Let $R=\bigcup_{n=1}^{\infty} \CC[[ t^{1/n} ]]$, the valuation ring of $K$. Then we can take the closure of $\overline{T x}$ in $G(d,n)(R)$ and take the fiber over $\Spec \CC$; geometrically this should be thought of as the limit of $\overline{T x}$ as $t \to 0$. Denote the fiber over $\Spec \CC$ by $Y$. $Y$ is itself a union of the toric varieties associated to the facets of $\DD$ -- see proposition \ref{ToricFiber}.

\begin{MT}
Suppose that $P$ is a tropical Pl\"ucker vector arising as $v(p_I(x))$ for some $x \in G(d,n)(K)$. Then  $\DD_P$ has at most $\frac{(n-c-1)!}{(d-c)! (n-d-c)! (c-1)!}$ interior faces of dimension $n-c$, with equality if and only if all of the facets of $\DD_p$ correspond to series-parallel matroids.
\end{MT}

Before explaining our strategy for proving the main theorem, we pause to describe some elegant results of Kapranov and of Hacking, Keel and Tevelev where $Y$ appears.  In \cite{CQoG1}, Kapranov defines a variety $X(d,n)$ called the Chow
quotient of the Grassmannian. The following description of $X(d,n)$ is equivalent to Kapranov's. Let $\Go(d,n)$ denote the open
subvariety of the Grassmannian where all $\binom{n}{d}$ Pl\"ucker
coordinates are nonzero. The torus $T$ acts on $G(d,n)$ and acts freely on
$\Go(d,n)$. $X(d,n)$ will be a certain compactification of $\Xo(d,n):=\Go(d,n)/T$.
We construct $X(d,n)$ as follows: for $x \in \Go(d,n)$, the closure
of the torus orbit through $x$, denoted $\overline{Tx}$, depends only
on the image of $x \in \Xo(d,n)$. Thus, $x \mapsto [\overline{Tx}]$
gives a map from $\Xo(d,n)$ to the Hilbert scheme of $G(d,n)$. We define $X(d,n)$
as the closure of the image of that map. (Kapranov uses the
Chow variety in place of the Hilbert scheme, but his Theorem 1.5.2
shows that this gives the same result.)  Kapranov shows that $X(2,n)$ is the moduli space $\overline{M_{0,n}}$ of stable genus zero curves with $n$ marked points and constructs maps between the $X(d,n)$ similar to the deletion maps between moduli spaces of stable curves.

By definition, each point $x$ of $X(d,n)$ corresponds to a $T$-invariant
subscheme of $G(d,n)$. If $x$ is a $K$-valued point of $\Xo(d,n)$, with preimage $(p_I) \in \Go(d,n)(K)$, then, as $X(d,n)$ is proper, we can extend $x$ to a map $\Spec R \to X(d,n)$. Let $x_0 \in X(d,n)(\CC)$ be the image of $\Spec \CC$ under this map and let $Y$ be formed as above from the Pl\"ucker coordinates $(p_I) \in K^{\binom{n}{d}}$. Then $Y$ is the $T$-invariant subscheme of $G(d,n)$ corresponding to $x_0$.

It is well known that $\Go(d,n)/T$ is the moduli space of arrangements of $n$ hyperplanes in $\PP^{d-1}$ in general position, considered up to automorphisms of $\PP^{d-1}$. This is sometimes called the Gelfand-MacPhearson correspondence, see \cite{GelfMac}. One might hope, then, that the points of $X(d,n)$ parameterize some sort of generalized hyperplane arrangements. In \cite{HKT}, Hacking, Keel and Tevelev explain how to do this: let $Y$ be the subscheme of $G(d,n)$ constructed above. For $e \in \CC^n \setminus \{ 0 \}$, let $G(d-1,n-1)_e$ be the subvariety of $G(d,n)$ corresponding to $d$-dimensional subspaces of $\CC^n$ containing $e$. Then $Z:=Y \cap G(d-1,n-1)_{(1,\ldots,1)}$ is a $(d-1)$-dimensional scheme and $Z \cap G(d-1, n-1)_{e_i}$, $i=1$, \dots $n$ provide $n$ hypersurfaces in  $Z$. Hacking, Keel and Tevelev show this is, in a sense motivated by log geometry, the correct generalization of an arrangement of $n$ generic hyperplanes in $\PP^{d-1}$. There is a stratification of  $Z$ by closed subvarieties in which the $d-i$ dimensional strata are in bijection with the $n-i$ dimensional interior faces of $\DD$. Our result is  thus a bound on the possible complexity of Hacking, Keel and Tevelev's ``very stable pairs".

We now turn to summarizing how we prove our result and how, in the process, we have discovered a matroid invariant which we hope to be of independent combinatorial interest. Let $Y$ be the subscheme of $G(d,n)$ described above. We mentioned above that $Y$ is a union of
torus orbit closures $Y=\bigcup_j \overline{T x^1_j}$, one for each facet of $\DD$. The
$\overline{T x_j^1}$ are glued along smaller torus orbit closures indexed by the smaller interior faces of $\DD$. More precisely, the structure sheaf $\OO_Y$ of $Y$ fits into an exact complex of sheaves
\begin{equation}
0 \ra \OO_Y \ra \bigoplus_{j=1}^{f_1} \OO_{\overline{T x_j^1}} \ra \cdots
\ra \bigoplus_{j=1}^{f_c} \OO_{\overline{T x_j^c}} \ra \cdots \ra
\bigoplus_{j=1}^{f_n} \OO_{\overline{T x_j^n} } \ra 0 \label{resolve}
\end{equation}
  where each $x_j^c$ is a point of $G(d,n)$ such that the torus orbit
  $T x_j^c$ is $n-c$ dimensional (see proposition \ref{ToricComplex}). We write $f_c$ for the number of strata of dimension $n-c$. 

In this language, our Main Theorem is
\begin{equation}
f_c \leq \frac{(n-c-1)!}{(d-c)!(n-d-c)!(c-1)!} \label{KeyEquation}
\end{equation}
with equality if and only if each $\overline{T x_i^1}$ corresponds to a series-parallel matroid.

The complex (\ref{resolve}) gives rise to an equality
\begin{equation}
[\OO_Y] = \sum_c (-1)^{c-1} \sum_{j=1}^{f_c} [\OO_{\overline{T x_j^c}}] \label{KTheoryEquality}
\end{equation}
in the $K$-theory $\Ko(G(d,n))$ of the Grassmannian. Moreover, $[\OO_Y]$ is the same as the class of $[\OO_{\overline{T x_{\mathrm{gen}}}}]$ where $x_{\mathrm{gen}}$ is a generic point of $G(d,n)$, because $Y$ is a flat degeneration of such a $\overline{T x_{\mathrm{gen}}}$. Thus, $[\OO_Y]$ depends only on $(d,n)$ and not on the particular decomposition $\DD$. Any additive linear map from $K$-theory to an abelian group gives restrictions on what components $\overline{T x_j^c}$ can occur in $Y$ and how many of them can occur.  However, because of the alternating signs in (\ref{KTheoryEquality}), most of these maps give no useful restrictions on the $f_c$. In the next two sections, we will construct a linear map $\Ko(G(d,n)) \to \ZZ[t]$ which has the positivity properties necessary to prove our Main Theorem. We denote the image of $[ \OO_{\overline{Tx}}]$ under this map by $h_x(t)$ and define $g_x(t):=(-1)^c h_x(-t)$ where $n-c=\dim Tx$. 

It turns out that $[\OO_{\overline{T x}}]$ depends only which Pl\"ucker coordinates of $x$ are nonzero -- data which is known as the matroid of $x$ (see proposition \ref{MatroidOnly}). Thus, $g_x$ is an invariant of realizable matroids and we will show in proposition \ref{NonRealizable} how to extend the definition to all matroids. We will thus, on occasion, fell free to write $g_M$ for $g_x$ where $M$ is the matroid whose bases are those $I \in \binom{[n]}{d}$ for which $p_I(x) \neq 0$.  $g_M$ turns out to be a very interesting invariant; it has
simple behavior under direct sum, two sum, duality, series and
parallel extension. Indeed, it is my hope that $g_M$ will emerge as a powerful invariant with many combinatorial descriptions.

In the next four paragraphs,  we describe one way of thinking about $g_M$ that combinatorialists may find particularly helpful. We emphasize that this is not a definition of $g_M$ but rather a description that can be deduced after carrying out all of the arguments in this paper.  For $M$ any matroid, let $\Poly_M$ denote the convex hull of the vectors $e_{i_1}+\cdots +e_{i_d}$ where $(i_1, \ldots, i_d)$ runs over the bases of $M$. Let $\DD$ be a  subdivsion of this polytope into subpolytopes which are also of the form $\Poly_N$ for various matroids $N$; we shall call such a subdivision \emph{matroidal}. Then we have $g_M=\sum g_N$ where the sum runs over $N$ such that $P_N$ is an interior face of $\DD$.  We will call a matroidal subdivision $\DD$ of $\Poly_M$ \emph{series-parallel} if every interior face of $\DD$ is a direct sum of series-parallel matroids. (It is equivalent to require that each facet of $\DD$ is series-parallel.) Corollary \ref{SerParGivesX} states that, whenever $N$ is a direct sum of $c$ series-parallel matroids, $g_N(t)=t^c$. We immediately obtain the following combinatorial corollary, which I believe is new.

\begin{prop}
Suppose that $\DD$ and $\DD'$ are two polyhedral subdivisions of $\Poly_M$, both of which have the property that all their interior faces correspond to direct sums of series-parallel matroids. Then $\DD$ and $\DD'$ have the same number of interior faces of each dimension.
\end{prop}

Note that there are many different series-parallel matroids for any given $(d,n)$. For example, there are 110 nonisomorphic series-parallel matroids of rank 5 on 10 elements. The polytopes $\Poly_M$ for the various series-parallel matroids of given $(d,n)$ do not have the same volume, same number of vertices, or have any other any other obvious invariant in common.

We can thus view $g_M=\sum_{i=1}^d g_i t^i$ as a generating function where $g_i$ is the number of faces of dimension $n-i$ occurring in any series-parallel decomposition of $\Poly_M$. Now, there are many matroids $M$ for which $\Poly_M$ can not be subdivided into series-parallel matroids -- the smallest example is the graphical matroid of the complete graph on four vertices. Nonetheless, we can consider $g_M$ as telling us how many series-parallel matroids $\Poly_M$ is morally made out of. It is already surprising that such a quantity can be consistently defined at all; it is further surprising that  (at least when $M$ is realizable over $\CC$) the coefficients of $g_M$ are nonnegative. This paper should be viewed as a challenge to combinatorialists --  what are the coefficients of $g_M$ counting when no series-parallel decomposition of $\Poly_M$ exists?

We now give some notation related to Grassmannians
and linear spaces. The standard basis of $\CC^n$ is called $e_1$, \dots, $e_n$. If $x \in G(d,n)$ then $L(x)$ denotes the
corresponding $d$-plane in $\CC^n$. Sometimes it will be convenient to
index the basis of $\CC^n$ by a finite set $A$ other than $[n]$; in
this case we write $G(d,A)$ for the Grassmannian, $\CC^A$ for the
vector space in which the $d$-planes live and $e_a$, $a \in A$, for the standard basis of $\CC^A$. $p_I(x)$ is the Pl\"ucker coordinate of $x$ indexed by $I=(i_1,\ldots,i_d) \in \binom{[n]}{d}$.
$M(x)$ is the matroid on $[n]$ whose bases are the $d$-tuples
$I=(i_1,\ldots,i_d)$ for which $p_I(x) \neq 0$. For $L \subset \CC^n$,
$L^{\perp}$ denotes the orthogonal complement of $L$ in $\CC^n$ under
the standard inner product $\langle (x_1,\ldots,x_n), (y_1,\ldots,y_n)
\rangle=\sum x_i y_i$. We abbreviate $L^{-1}(L(x)^{\perp})$ by
$x^{\perp}$. If $v \in \CC^n \setminus \{ 0 \}$, then $G(d-1,n-1)_v$
denotes the subscheme of $G(d,n)$ consisting of those $x$ for which $v
\in L(x)$ and $G(d,n-1)^v$ consists of those $x$ for which $v \in
L(x)^{\perp}$. We define the $2n$ subschemes of the form
$G(d-1,n-1)_{e_{i}}$ and $G(d,n-1)^{e_i}$ to be the ``coordinate
subgrassmannians'' of $G(d,n)$. Many of our theorems will have as a
hypothesis that some point not lie in any coordinate subgrassmannian.
It is easy to check that the $x^c_j$ occuring in the complex (\ref{resolve}) do not lie in any coordinate subgrassmannian.

At this point, let me give some reassurance to algebraic geometers not familiar with matroid theory. For $x \in G(d,n)$, the matroid associated to $x$ is the ordered pair $( [n], B )$ where $B$ is the set of $I \in \binom{[n]}{d}$ for which $p_I(x) \neq 0$.  (See chapter 2 of \cite{White1} for details and other, equivalent, definitions.) A matroid in general is an ordered pair $(S,B)$ where $B$ is a subset of $\binom{S}{d}$ which obeys some combinatorial axioms; these axioms are easily checked in the case of a matroid which comes from a point of $G(d,n)$. Matroid theorists have developed an extremely useful vocabulary for describing operations under which we combine matroids to produce new matroids. All such operations which we will use in this paper correspond to actual geometric operations on Grassmannians. Our practice will always be to give a definition of the operation on Grassmannians and use notation for it as similar as possible to the notation used by matroid theorists. Thus, the reader not familiar with matroid theory should always be able to follow the main flow of the paper by simply using our geometric definition. We will, however, feel free to use matroid terminology in offhand comments before giving the corresponding definition.

The ideas for this paper began to ferment when working on my dissertation under the supervision of Bernd Sturmfels. His guidance has greatly helped me understand $K$-theory, toric varieties and matroids and has given me the concrete tools necessary to compute with them. I have also benefitted from conversations with Paul Hacking, Allen Knutson, Ezra Miller and Vic Reiner and I might never have finished the proof of proposition \ref{Positivity} if  Robert Lazarsfeld had not explained to me how to think about Kawamata-Viehweg vanishing. Finally, I am glad to acknowledge the financial support of the Clay Foundation during the research and writing of this paper.

\section{Geometric Preliminaries}

In this section, we will review the material on $K$-theory which we will need for this paper. Our strategy is strongly influenced by the methods of \cite{Brion} and that paper will also serve a good reference for the material in this section.

Let $X$ be a smooth variety over an algebraically closed field $k$. Let $\Ko(X)$ denote the Grothendieck group of coherent
sheaves on $X$ -- because we will only work with smooth varieties, this is the same as the Grothendieck group of vector bundles. If $\EE$ is a coherent sheaf on $X$, we
denote its class in $\Ko(X)$ by $[\EE]$. If $f:X_1 \to X_2$ is a flat map of smooth varieties we get a pull back map $f^* \Ko(X_2) \to \Ko(X_1)$ given by $f^* [\EE] =[f^* \EE]$. If $f : X_1 \into X_2$ is a closed inclusion, we again get a pull back map given by $f^* : [\EE] = \sum_i (-1)^i [\TorTor_i(\EE, \OO_{f(X_1)})]$. If $f$ can be written as $\pi \circ \iota$ for a flat map $\pi$ and a closed inclusion $\iota$ (\emph{e.g.} if $f$ is projective), one can check that $\iota^* \pi^*$ depends only on $f$, and we denote this map as $f^*$. In fact, $\TorTor_i(\pi^* \EE, \OO_{\iota(X_1)})$ depends only on $f$ and will be denoted $\TorTor_i^{f}(\EE,\OO_{X_1})$.

If $f$ is a proper map, we also get a pushforward $f_* : \Ko(X_1) \to \Ko(X_2)$ given by $f_* [\EE] = \sum_{i} (-1)^i [R^i f_*(\EE)]$. If $f$ is finite this is simply $[\EE] \mapsto [f_* \EE]$. The additive group $\Ko(X)$ is made into a ring with multiplication
$$[\EE] [\FF] = \sum_{i} (-1)^i [\TorTor_i(\EE, \FF)].$$
The maps $f^*$ and $f_*$ are both functorial and additive, but only $f^*$ is a map of rings.

If $X$ is proper then pushforward gives us a map $\Ko(X) \to \Ko(\Spec k)=\ZZ$. This map is simply
$[\EE] \mapsto \sum (-1)^i \dim_k H^i(\EE, X)$ and is denoted by $\chi$.  When $\EE$ is the structure sheaf of a closed subscheme $Z \subset X$, we will abbreviate $\chi(\OO_Z)$ by
$\chi(Z)$. Note that $\chi(Z)$ depends only on $Z$ as an abstract
scheme and not on its embedding into $X$. In this setting, we
will call $\chi(Z)$ the holomorphic Euler characteristic of $Z$. (The word ``holomorphic" is inserted to
distinguish $\chi(Z)$ from the topological Euler characteristic of
$Z(\CC)$.)

In this paper, most of the varieties we will consider have large symmetry groups. In that context, the following result is useful.

\begin{prop} \label{BertiniTor}
Suppose that $G$ is a connected linear algebraic group acting transitively on $X$. Then $G$ acts trivially on $\Ko(X)$. Given coherent sheaves $\EE$ and $\FF$ on $X$, for a generic $g \in G$ we have $\TorTor_i(\EE, g \FF)=0$ for $i>0$. Similarly, let $f:X_1 \to X$ be a projective map and let $\EE$ be a coherent sheaf on $X$. Then, for a generic $g \in G$, we have $\TorTor_i^{f}(g \EE, \OO_{X_1})=0$.
\end{prop}

\begin{proof} 
Let $\GG$ be the coherent sheaf on $G \times X$ obtained by pulling back $\EE$ along the multiplication map $G \times X \to X$. Then $[g \EE]$ is $(\pi_2)_* ([\GG] \pi_1^*[k_g])$ where $\pi_1$ and $\pi_2$ are the projections of $G \times X$ onto its components and $k_g$ is the skyscraper sheaf at $g \in G$. For any $g_1$ and $g_2 \in G$, there is a rational curve $C \subset G$ containing $g_1$ and $g_2$.  The skyscraper sheaves at $g_1$ and $g_2$ are equivalent in $\Ko(C)$ so, pushing forward this equality, we get $[k_{g_1}]=[k_{g_2}]$ in $\Ko(G)$ and we get $[g_1 \EE]=[g_2 \EE]$. We have now shown that $G$ acts trivially on $\Ko(X)$.

The claim about $\TorTor_i(\EE, g \FF)$ vanishing is the main result of \cite{MS}. To extend this to the second $\TorTor$ vanishing result, let $f = \pi \circ \iota$ where $\iota : X_1 \into X \times \PP^N$ is a closed immersion and $\pi : X \times \PP^N \to X$ is the projection. Then $G \times \PGL_{N+1}$ acts transitively on $X \times \PP^N$ and the $\TorTor$ we are trying to compute is $\TorTor_i(\iota_*(\OO_{X_1}), (g,h) \pi^* \EE)$ where $h \in \PGL_{N+1}$ is chosen arbitrarily. So our second $\TorTor$ vanishing claim is reduced to the first.
\end{proof}

Let $Y$ and  $x_j^c$ be as in the preceeding
section. Let $x_{\gen}$ be some point of 
$\Go(d,n)$; the class $[\OO_{\overline{T x_{\gen}}}] \in \Ko(G(d,n))$ does not depend on the choice of $x_{\gen}$. Then $Y$ is a flat degeneration of $\overline{T x_{\gen}}$, $[\OO_Y]=[\OO_{\overline{ T x_{\gen}}}]$ and, by equation (\ref{KTheoryEquality}), we have
$$[\OO_{\overline{T x_{\gen}}}]=\sum_c (-1)^{c-1} \sum_j [\OO_{\overline{T x_j^c}}].$$
Our
method of proving the Main Theorem will be to pair both sides of this
equality with certain Schubert classes $\Omega_{\lambda}$. By the above observations, for
$\lambda$ any partition fitting inside a $d \times (n-d)$ box and $g$
a generic member of $\PGL_n$,
$$\sum_c (-1)^{c-1} \sum_j \chi(g \OO_{\Omega_{\lambda}} \otimes
\OO_{\overline{T x_j^c}}) \label{key}$$
is a constant depending only
on $\lambda$, $d$ and $n$, not on the particular degnerate fiber $Y$. Now, $\OO_{g \Omega_{\lambda}} \otimes \OO_{\overline{T
    x_j^i}} = \OO_{g \Omega_{\lambda} \cap \overline{T x_j^i}}$ where
the intersection is taken in the scheme theoretic sense. 

We will show in proposition \ref{RationalSing} that $\overline{Tx}$ has rational singularities. The reason that this is important is the following result:

\begin{prop} 
Let $Y$ and $Y'$ be proper varieties with rational singularities which are birational to each other. Then $\chi(\OO_{Y})=\chi(\OO_{Y'})$.
\end{prop}

\begin{proof}
  Let $\pi: \tilde{Y} \to Y$ and $\pi' : \tilde{Y'} \to Y'$ be
  desingularizations of $Y$ and $Y'$ respectively. The definition of
  rational singularities is that $\pi_* \OO_{\tilde{Y}}=\OO_Y$ and $R^i
  \pi_* \OO_{\tilde{Y}}=0$. For any coherent sheaf $\EE$ on $\tilde{Y}$,
  there is a spectral sequence relating $H^j(R^i \pi_* \EE)$ to
  $H^k(\EE)$; in this case, we get that
  $H^i(\OO_{\tilde{Y}},\tilde{Y})=H^i(\OO_Y,Y)$ and so
  $\chi(\tilde{Y})=\chi(Y)$ and, in the same way, $\chi(\tilde{Y'})=\chi(Y')$. We
  are thus reduced to showing $\chi(\tilde{Y})=\chi(\tilde{Y'})$, that
  is, proving the same result for $Y$ and $Y'$ smooth. This is a
  classical result.
\end{proof}

\begin{cor} \label{Birational}
  Let $g$ be a generic element of $\PGL_n(\CC)$. Let $W$ be a proper smooth
  variety birational to $\overline{T x} \cap g \Omega_{\lambda}$. Then
  $\chi(W)=\chi(\overline{T x} \cap g \Omega_{\lambda})$.
\end{cor}

\begin{proof}
  By proposition \ref{RationalSing}, $\overline{Tx}$ has rational singularities.  All Schubert varieties have rational singularities
  so, by lemma 2 of \cite{Brion}, $\overline{T x} \cap g
  \Omega_{\lambda}$ has rational singularities for generic $g$. Then
  the above proposition tells us that $\chi(W)=\chi(\overline{T x}
  \cap g \Omega_{\lambda})$.
\end{proof}


\section{Construction of $g_x$ and summary of results}

We now describe the construction of the map $K^\circ(G(d,n)) \to \ZZ[t]$ that we promised in section \ref{Intro}. Let $i$ be
an integer between $1$ and $n-1$. Choose a generic \footnote{We will frequently use the adjective generic in this paper. In each case it means that the results that follow are true as long as the object in question is chosen in some Zariski open subset of the obvious parameter space of objects.} $n-i$
plane $M_i \subset \CC^n$ and a generic line $\ell \subset M_i$; we
define $\Omega_i \subset G(d,n)$ to be the Schubert variety consisting of those $x$
such that $\ell \subset L(x)$ and $L(x) + M_i \neq \CC^n$. Note that this second condition is vacuous when $i > d$, so $\Omega_{d+1}=\Omega_{d+2}=\cdots=\Omega_{n-1}$. We  extend our notation by setting $\Omega_{N}=\Omega_{d+1}$ when $N \geq n$.  Let $Z$ be a
closed subvariety of $G(d,n)$. We define a formal power series $r_Z$
by
$$r_Z(t)=\sum_{i=1}^{\infty} \chi([\OO_Z] [\Omega_i])
t^i.$$
Since the coefficient of $t^i$ becomes constant for $i$
sufficiently large, $r_Z(t)$ is a rational function of the form
$h_Z(t)/(1-t)$. Applying this operator to both sides of equation
(\ref{KTheoryEquality}), we get
$$h_{Y}(t)=\sum_{c=1}^n (-1)^{c-1} \sum_{j=1}^{f_i} h_{\overline{T x_j^i}}(t). \label{key2}$$
For $x \in G(d,n)$, $x$ not in any coordinate subgrassmannian, define 
$$g_x(t)=(-1)^c h_{\overline{T x}}(-t)$$
where $n-c = \dim T x$.

\textbf{Remark:} Generally speaking, combinatorial results are slightly nicer when stated in terms of $g$ and geometric results are slightly nicer when stated in terms of $h$. In this paper, we favor $g$.

The Main Theorem will follow from several lemmas about the behavior of $g_x$ which we now list.

\begin{prop} \label{GenericCase}
For $x_{\textrm{gen}} \in G(d,n)$ a point none of whose Pl\"ucker coordinates are zero, we have
$$g_{x_{\textrm{gen}}}(t)=\sum_{i=1}^d \frac{(n-i-1)!}{(d-i)!(n-d-i)!(i-1)!} t^i.$$
\end{prop}

Let $x_1 \in G(d_1,n_1)$, \dots, $x_r \in G(d_r,n_r)$. Define $\bigoplus x_k \in
G(\sum d_k, \sum n_k)$ to be the point corresponding to the direct sum
$\bigoplus L(x_k) \subset \CC^{\sum n_k}$. 

\begin{prop} \label{ProductFormula}
With notation as above, we have
$$g_{\bigoplus x_k}(t) = \prod g_{x_k}(t).$$
\end{prop}

\begin{prop} \label{Positivity}
Assume that $x$ is not contained in any coordinate subgrassmannian and $n \geq 2$. Then the coefficents of $g_x$ are nonnegative and the coefficient of $t^c$ in $g_x$ is positive, where $n-c=\dim Tx$.
\end{prop}

\textbf{Remark:} In \cite{Brion}, Brion proves that certain linear combinations of the quantities $\chi([W][\Omega_{\lambda}])$ are nonnegative for any subscheme $W \subset G(d,n)$ with rational singularities. These inequalities are not strong enough to imply proposition \ref{Positivity}.

The proofs of these propositions will occupy sections \ref{BetaInvariant} through
\ref{PosProof}. For now, let us see why they imply the Main Theorem. From
equation (\ref{key2}) and proposition \ref{GenericCase}, we have
$$\sum_{i=1}^d (-1)^{i+1} \frac{(n-i-1)!}{(d-i)!(n-d-i)!(i-1)!}  t^i = \sum_c (-1)^{c+1} \sum_{j=1}^{f_c} h_{\overline{T x_j^c}}(t).$$
Substituting $-t$ for $t$ and negating both sides, we get
$$\sum_{i=1}^d \frac{(n-i-1)!}{(d-i)!(n-d-i)!(i-1)!} t^i = \sum_c  \sum_{j=1}^{f_c} g_{x_j^c}(t).$$
By proposition \ref{Positivity}, every term of the polynomials on the right hand
side is nonnegative and the $t^c$ term of $g_{x_j^c}(t)$ is positive, so
$$\sum_{i=1}^d \frac{(n-i-1)!}{(d-i)!(n-d-i)!(i-1)!} t^i \succeq
\sum_c \sum_{j=1}^{f_c} t^c=\sum_c f_c t^c$$
where $\succeq$ denotes term by term dominance.
In other words, $\frac{(n-i-1)!}{(d-i)!(n-d-i)!(i-1)!} \geq f_i$, exactly as we wanted. \qedsymbol

Note that, in particular, we have shown that $0 \geq f_i$ when $i > \min(d,n-d)$, so we have shown that complex (\ref{resolve}) stops after $\min(d,n-d)$ steps.

We now summarize the rest of the paper. In section \ref{DecompositionResult}, we prove a technical result that will allow us to reduce many of our arguments to the case $c=1$. In section \ref{BetaInvariant},  we prove the positivity of the leading term of $g_x$; this result is not only a special case of proposition \ref{Positivity} but is used in a crucial way in the proof of that proposition. In section \ref{GammaConstruction}, we introduce a geometric construction that is extremely useful for proving results about $g_x$. In addition, we show that $g_x$ is invariant under orthogonal complement, and under series and parallel extensions. In sections \ref{ProdProof} and \ref{PosProof}, we prove propositions \ref{ProductFormula} and \ref{Positivity}. At this point, we will have essentially proven the Main Theorem.  (We delay the proof of proposition \ref{GenericCase} to section \ref{Examples}, but this is only a computation.) We now switch to the question of computing $g_x$. In section \ref{SectionTwoSum}, we show that $g_{x_1 +_2 x_2}(t)=g_{x_1}(t) g_{x_2}(t)/t$. This allows us to reduce the computation of $g_x$ in many cases to the computation of $g_{x'}$ for simpler $x'$. In section \ref{Examples}, we give many examples in which we compute $g_x$. We close with numerous speculations and conjectures. At the end of the paper, we have included an appendix which proves some basic facs about torus orbits in $G(d,n)$. 

\section{A decomposition result} \label{DecompositionResult}

In this section we will prove a result that will let us reduce many of our results to the
case $c=1$.  This result is widely known, but it is usually stated in
the language of matroids so one must then unwrap the matroid
definitions. It seems simplest to give a proof.

\begin{prop} \label{Decomposition}
  Let $\overline{T x}$ have dimension $n-c$ and assume that $x$ is not
  contained in any coordinate subgrassmannian. Then there
  exists a partition $[n]=\bigsqcup_{k=1}^c A_k$ of $[n]$ into $c$ parts, positive
  integers $d_1$, \ldots, $d_c$ with $\sum d_k=d$ and $d_k < |A_k|$ and
  points $x(1) \in G(d_1,A_1)$, \dots, $x(c) \in G(d_c,A_c)$ such that
  $\dim \overline{T(k) x(k)}=|A_k|-1$ and $x=\bigoplus x(k)$.
\end{prop}

Here $T(k)$ denotes the torus $(\CC^*)^{A_k}/\diag(\CC^*)$.

\begin{proof}
Recall that $I \in \binom{[n]}{d}$ is called a basis of $M(x)$ if the
Pl\"ucker coordinate $p_I(x)$ is nonzero. The (complex) dimension of
$\overline{T x}$ is the same as the (real) dimension of its moment map image $P_{M(x)}$.
The polytope $P_{M(x)}$ is the convex hull of the vectors $e_{i_1}+\cdots +
e_{i_d}$ where $(i_1,\ldots,i_d)$ ranges over the bases of $M(x)$. As was observed by Gelfand,
Goresky, MacPhearson and Serganova \cite{GGMS}, all of the edges of
$P_{M(x)}$ are parallel to $e_i-e_j$ for some $1 \leq i < j \leq
n$. 

The dimension of $P_{M(x)}$ is the same as the dimension of the affine
linear space $L$ it spans which is, in turn, the same as the dimension of the vector
space $V$ generated by the directions of the edges of $P_{M(x)}$. So we
must compute the dimension of a vector space spanned by vectors of the
form $e_i-e_j$. Define an equivalence relation $\sym$ on $[n]$ to be
generated by the relations $i \sym j$ if there is an edge of
$P_{M(x)}$ parallel to $e_i-e_j$. Let $A_1$, \dots, $A_s$ be the
equivalence classes of $\sym$. Then $V$ has dimension $n-s$, so $s=c$, and $V$ is cut
out by the equations $\sum_{j \in A_r} x_j=0$ for $1 \leq r \leq c$.
$L$ is cut out by equations of the form $\sum_{j \in A_r}
x_j=\textrm{constant}$. Take this constant to be $d_r$.

We clearly have $\sum d_r = d$ and $\bigsqcup A_r=[n]$. Write $\CC^n =
\bigoplus \CC^{A_r}$ and let $L(x_r)=L(x) \cap \CC^{A_r}$. We claim
that $L(x_r)$ is $d_r$-dimensional. Proof: Every basis of $M(x)$ contains
exactly, and in particular no more than, $d-d_r$ elements of $[n] \setminus A_r$.  This implies that the
dimension of the projection of $L(x)$ to $\CC^{[n] \setminus A_r}$ is
at most $d-d_r$. Thus $\dim L(x_r) \geq d_r$. But $\sum d_r=d$ and,
since $\CC^n = \bigoplus \CC^{A_r}$, we must have $\sum \dim L(x_r)
\leq \dim L(x)=d$. So we have equality and $\dim L(x_r)=d_r$. In order
to have equality, we must have $L(x)=\bigoplus L(x_r)$.

The strict inequalities in $0 < d_r < |A_r|$ follow from the
assumption that $x$ is not contained in a coordinate subgrassmannian.
\end{proof}

\textbf{Remark:} The sets $A_i$, equipped with the structure of a matroid by the points $x_i \in G(d_i,A_i)$ are called the connected components of $M(x)$. We will always use $c$ to denote the number of connected components of $M(x)$. $M$ is called connected if $c=1$. See section 6.2 of \cite{White1} for more on the connected components of matroids.

\section{The $\beta$-invariant and cohomology} \label{BetaInvariant}

In this section, we will compute $\chi(\OO_{\overline{T x} \cap
  \Omega_1})$. By considerations of dimension, we see that
$\overline{T x} \cap \Omega_1$ is empty when $\dim \overline{T x} <
n-1$ and finite when $\dim \overline{T x} = n-1$. Thus, $\chi$ simply
counts the number of points of $\overline{T x} \cap \Omega_1$. In
other words, we are being asked to determine, given a generic
hyperplane $H$ and a generic line $\ell \subseteq H$, for how many
points $y \in \overline{T x}$ we have $\ell \subset L(y) \subset
H$.  We will denote the set of $y$ such that $\ell \subset L(y)
\subset H$ by $\Omega(\ell,H)$. 

The result we will establish in this section is that the size of
$\overline{T x} \cap \Omega(\ell,H)$ is equal to a well known
combinatorial invariant, the $\beta$ or Crapo invariant, of the
matroid $M(x)$. $\beta(M)$ is one of the best known invariants of a
matroid $M$, see chapter 6 of \cite{White2} for a survey of its significance. One can define $\beta$ by the fact that it obeys the Tutte recurrence
$\beta(M)=\beta(M/e)+\beta(M \setminus e)$ for $|M| \geq 3$, that
$\beta(M)=0$ if $M$ has a loop or coloop and $\beta(M)=1$ if $M$ is
the uniform matroid of rank $1$ on $2$ elements. It is not clear that there is a well defined matroid invariant with these properties; one may consider this section to be a geometric proof.

\begin{theorem} \label{Beta}
  Let $H \subset \CC^n$ be a generic hyperplane and $\ell \in H$ a
  generic line in $H$.  Then $\#(\overline{T x} \cap \Omega(\ell,H))$ is $\beta(M(x))$.  
\end{theorem}

For $x \in G(d,n)$, $n \geq 2$, we denote the value of
$\#(\overline{T x} \cap \Omega(\ell,H))$ for generic $(\ell,H)$ by $b(x)$. (When $n=1$, this
formula doesn't make sense as $\dim \ell =1 > \dim H=0$ so it is
impossible to find $\ell \subset H$.) It is enough to show that this number obeys the defining recurrences
of the $\beta$-invariant. We now cast each of these into a
geometric statement and prove it.

We introduce the following notation: let $x \in G(d,n)$ and let $i \in
[n]$. Assuming that $L(x) \not \subseteq \{ z_i=0 \}$, we define $x/i \in G(d-1,[n] \setminus \{ i \})$ so that
$L(x/i)=L(x) \cap \{ z_i=0\}$. Assuming that $L(x)$ does not contain
$e_i$, we define $x \setminus i \in G(d,[n] \setminus \{ i \})$ so that
$L(x \setminus i)$ is the image of $L(x)$ in $\CC^n/e_i$.

\begin{prop} \label{Additive}
Let $n \geq 3$ and $i \in [n]$. Let $x \in G(d,n)$ and assume that $L(x) \not \subseteq \{ z_i=0 \}$ and $e_i \not \in L(x)$. Then $b(x)=b(x/i)+b(x \setminus i)$.
\end{prop}

\begin{proof}
  Let $\ell$ be a generic line in the hyperplane $z_i=0$ and let $H'$ be a
  generic hyperplane in $z_i=0$ containing $\ell$. Let $H=H'
  \oplus e_i$. Of course, $(\ell,H)$ is not a generic pair ``line,
  hyperplane containing line'' in $\CC^n$. Nonetheless, we claim that
  $\overline{Tx}$ meets $\Omega(\ell,H)$ transversely and that it does
  so at $b(x/i)+b(x \setminus i)$ points. 
  
  We divide $\overline{Tx}$ into three pieces: a closed piece
  $X_1$ consisting of those $y \in \overline{T x}$ for which $L(y)
  \subseteq \{z_i=0 \}$, a closed piece $X_2$ consisting of those $y \in
  \overline{T x}$ for which $L(y) \ni e_i$ and an open piece $U$ which
  is the complement of $X_1 \cup X_2$. Note that $X_1 \cap
  X_2=\emptyset$. We claim that $\Omega(\ell,H) \cap U=\emptyset$.
  Suppose on the contrary that $y \in \Omega(\ell,H) \cap U$. Then
  $y/i$ and $y \setminus i$ are both well defined. Consider $(y/i, y
  \setminus i)$ as a point of the two-step flag manifold of pairs ``$(d-1)$-plane,
  $d$-plane containing $(d-1)$-plane'' in $(n-1)$-space. Then $(y/i, y \setminus i)$ lies in the closure of the $T':=(\CC^*)^{[n] \setminus \{ i \}}/\diag \CC^*$ orbit through $(x/i,
  x \setminus i)$. Moreover, the condition that $y \in
  \Omega(\ell,H)$ is equivalent to ``$\ell \subset L(y/i)$ and $L(y
  \setminus i) \subset H'$''. Thus, to show our claim, we must show that
  for $\ell$ and $H'$ chosen generically according to the given
  constraints, there is no $(y_1,y_2)$  in $\overline{T' (x_1,x_2)}$ for which $\ell \subset L(y_1)$ and $L(y_2) \subset H'$. The $T'$-orbit
  closure is at most $(n-2)$-dimensional and, for each point $(y_1,y_2)$ in
  the $T'$-orbit closure, the space of possible choices for $\ell$ and
  $H'$ are $d-2$ and $(n-d-2)$-dimensional, respectively. Thus, we
  have at most $(n-2)+(d-2)+(n-d-2)=2n-6$ dimensions of possible pairs
  $(\ell, H')$ which are compatible with some $(y_1,y_2) \in \overline{T' (x/i, x \setminus i)}$. There are, on the other hand, $(n-2)+(n-3)=2n-5$
  dimensions of possible pairs $(\ell,H')$ with $\ell \subset H'
  \subset \CC^{n-1}$. Thus, for a generic $(\ell,H')$, $\Omega(\ell,H)
  \cap U=\emptyset$ as claimed.

We see that, as point sets, 
$$\Omega(\ell,H) \cap \overline{T x}=\left( \Omega(\ell,H) \cap X_1
\right) \sqcup \left( \Omega(\ell,H) \cap X_2 \right).$$
We will now show that
this is in fact true as an equality of schemes. This is a local
question, we check it on each of two open sets $V_1:=\{ y : L(y) \not \ni e_i \}$
and $V_2:=\{ y: L(y) \not \subseteq \{ x_i =0 \} \}$, which together form a cover of $G(d,n)$.
The variety $V_1$ is a $d$-dimensional vector bundle over
$G(d,n-1)$. The closed subschemes $X_1 \cup U$ and $\Omega(\ell,H) \cap V_1$ are each
sub-vector bundles over subschemes of $G(d,n-1)$ -- specifically, over
$X_1$ and $\Omega(\ell,H')$ respectively. So their intersection has a map to
$G(d,n-1)$ where each scheme theoretic fiber is a vector space. But we
know that this intersection is disjoint from $U$, which can only
happen if each of those fibers are zero-dimensional. This, in turn,
shows that the intersection is contained in $X_1$, not only on the
level of point sets, but scheme-theoretically. An analogous argument
show that the part of $\Omega_{\ell,H} \cap \overline{T x}$ in the
open set $V_2$ is scheme-theoretically contained in $X_2$.

Now, $X_1 = \overline{T' (x \setminus i)}$ and $\Omega(\ell,H) \cap
X_1=\Omega(\ell,H') \cap \overline{T' (x \setminus i)}$. So, by induction, for
generic $(\ell,H')$, the intersection $\Omega(\ell,H) \cap X_1$ is transverse and consists of $b(x
\setminus i)$ isolated points. Similarly, $X_2 \cap \Omega(\ell,H)$ is transverse and consists of $b(x/i)$
isolated points. In conclusion, the intersection $\Omega(\ell,H) \cap
\overline{Tx}$ is transverse and consists of $b(x \setminus
i)+b(x/i)$ points.

\end{proof}

\begin{prop} \label{Zero}
  If $L(x)$ is contained in one of the hyperplanes of the form $\CC^{[n] \setminus \{ i \}}$
  or contains one of the vectors $e_i$ then $b(x)=0$.
\end{prop}

\begin{proof}
  If $L(x) \subset \{ z_i=0 \}$ then $L(y) \subset \{ z_i=0 \}$ for
  all $y \in \overline{T x}$. As a generic line $\ell$ is not
  contained in $\{ z_i =0 \}$, $\overline{T x} \cap \Omega(\ell,H) =
  \emptyset$ for a generic $\ell$. Similarly, if $e_i \in L(x)$ then  $\overline{T x} \cap \Omega(\ell,H) =
  \emptyset$ for a generic $H$.
\end{proof}

\begin{prop} \label{Normalize}
If $x$ is a generic point of $G(1,2)$ then $b(x)=1$.
\end{prop}

\begin{proof}
  The Grassmannian $G(1,2)$ is just the projective line. Assuming $x$
  is a generic point on this line, $\overline{T x}=\PP^1$. We have
  $\ell=H$ and $\Omega(\ell,H)$ is simply a point. Thus, $\overline{Tx}
  \cap \Omega(\ell,H)$ is a single point as desired.
\end{proof}

Propositions \ref{Additive}, \ref{Zero} and \ref{Normalize} show that
$b(x)$ obeys the defining recurrence and initial conditions of the
$\beta$ invariant. Thus, theorem \ref{Beta} is proven. \qedsymbol

The key importance of this result for us will be that we can use it to
show that $\overline{T x} \cap \Omega_{1}$ is nonempty whenever $\dim
Tx$ is large enough. Specifically,

\begin{prop} \label{CoeffOfT}
  If $\overline{T x}$ is $(n-1)$-dimensional and $n \geq 2$ then $\overline{T x} \cap
  \Omega_{1}$ is nonempty and finite. As a corollary, the coefficient of $t$ in $g_x(t)$ is nonzero in this case.
  \end{prop}

\begin{proof}
  The second claim follows from the first, as $c_1=\chi(\overline{T x}
  \cap \Omega_1)$. For the first, we appeal to theorem II of \cite{Crapo}: if $M$ is a connected matroid with $n \geq 2$ then $\beta(M)>0$.
  That $M$ is connected precisely means that all of $[n]$ is a single
  equivalence class under the equivalence relation in the proof of
  proposition \ref{Decomposition}, which, by the proof of that
  proposition, is equivalent to saying that $\overline{Tx}$ is $n-1$
  dimensional.
\end{proof}

This result is not only important for establishing the positivity of
$c_1$; it will also be used to establish the generic finiteness of a
map in section \ref{PosProof}, which will in turn be used to allow us to
prove proposition \ref{Positivity} as a corollary of Kawamata-Viehweg vanishing.

It is difficult to give an attribution for theorem \ref{Beta} which is
why we have included a complete proof. At the same time, this result is not truly original. The following paragraphs explain
how theorem \ref{Beta} could be pieced together from previously published results. The
problem of computing $b(x)$ is related to the following problem from
algebraic statistics:

\textbf{Problem:} Let $a_1$, \ldots, $a_{n-1}$ be $n-1$ affine linear functionals on
$\CC^{d-1}$ with $\sum a_i=1$. Let $p_1$, \dots, $p_{n-1}$ be positive integers which
are generic (meaning that they are in the complement of the zero locus
of finitely many polynomials, this collection of polynomials depending
on the $a$'s). Compute the number of critical points of
$$\Phi(u):=\prod_{i=1}^{n-1} a_i(u)^{p_i} \label{MaxThis}$$
on $\CC^{d-1} \setminus \bigcup_{i=1}^{n-1} \{ a_i(u)=0 \}$. 

This problem arises naturally when there is some experiment whose
outcome depends on $d-1$ parameters $(u_1, \ldots, u_{d-1})$ with
unknown values and which can yield $n-1$ outcomes. If the probability
of outcome $i$ is $a_i(u)$ and
$p_i$ is the number of times that outcome $i$ was observed, it is
standard to estimate $u$ by maximizing $\Phi$ over real values of $u$.
Let $x \in G(d,n)$ be such that $L(x)$ is the $d$-plane in $n$-space with $L(x) \cap \{ z_n=1 \}$
parameterized by $(u_1, \ldots u_{d-1}) \mapsto (a_1(u) , \ldots,
a_{n-1}(u) , 1)$.  Then it can be shown that the critical points
in question are in bijection with the points of
$\overline{T x} \cap \Omega(\ell,H)$ where $\ell=\Span(1,\ldots,1)$
and $H=(p_1, \ldots, p_{n-1}, -\sum_{i=1}^{n-1} p_i)^{\perp}$.  I am
not aware of a reference which points out the connection between
this problem and the intersection theory problem of enumerating
$\overline{Tx} \cap \Omega(\ell,H)$.  However, theorem 28 of
\cite{MaxLikeDeg} describes the critical points of $\Phi$ as the top
chern class of a certain sheaf of logarithmic differentials on any
compactification $X$ of $U:=\CC^{d-1} \setminus \bigcup_{i=1}^{n-1} \{
a_i(u)=0 \}$ in which $X \setminus U$ becomes a normal crossing
divisor obeying certain conditions.  It is observed in section 2.2 of
\cite{HKT} that $\overline{T x} \cap G(d-1,n-1)_{\ell}$ is such a
compactification and that the sheaf of logarithmic differentials
involved is the restriction of the anti-tautological bundle of
$G(d-1,n-1)_{\ell}$. Since chern classes are contravariant, the
equivalence of the chern class description and the intersection theory
description is simply the standard fact that $\Omega(\ell,H) \subset
G(d-1,n-1)_{\ell}$ represents the top chern class of the
anti-tautological bundle. I am grateful to Paul Hacking for pointing out to me the connections described in this paragraph.

Varchenko \cite{Var} considered the problem of determining the number
of critical points of $\Phi$ and showed that, when $x \in
G(d,n)(\RR)$, the number of critical points of $\Phi$ is equal to the
number of bounded regions of $\RR^{d-1} \setminus \bigcup_{i=1}^{n-1}
\{ a_i(u)=0 \}$ -- in fact, there is exactly one critical point in each bounded region. The number of such regions is equal to
$\beta(M(x))$.  Varchenko also considered the case where $x$ is not
defined over $\RR$ and conjectured that in this case the number of
critical points is equal to the (topological) Euler characteristic of
$\CC^{d-1} \setminus \bigcup_{i=1}^{n-1} \{ a_i(u)=0 \}$; this Euler
characteristic is also known to be equal to $\beta(M(x))$.
Orlik and Terao proved Varchenko's conjecture correct in \cite{OT}.

\section{A Variety Birational to $\overline{Tx} \cap \Omega_i$} \label{GammaConstruction}

We saw in proposition \ref{Birational} that, before computing the holomorphic Euler characteristic of $\overline{Tx} \cap \Omega_i$, we may replace
$\overline{Tx } \cap \Omega_i$ by any smooth proper variety birational to it. In this section, we will present such a variety which
will be very useful for proving later results. 

Let $x \in G(d,n)$ and let $n-c=\dim \overline{Tx}$. We have $\PP(L(x)) \times
\PP(L(x)^{\perp}) \subset \PP^{n-1} \times \PP^{n-1}$. There is a
birational map $m: \PP^{n-1} \times \PP^{n-1} \dashrightarrow
\PP^{n-1}$ given by
$$m : (x_1 : \cdots : x_n) \times (y_1: \cdots : y_n) \mapsto (x_1 y_1
: \cdots : x_n y_n).$$
As long as $x$ is not contained in any coordinate subgrassmanian, a generic point of $\PP(L(x))
\times \PP(L(x)^{\perp})$ has all coordinates non-zero, so $m$ is
defined on a dense open subset of $\PP(L(x)) \times \PP(L(x)^{\perp})$.
Since $L(x)$ and $L(x)^{\perp}$ are perpendicular, $m\left( \PP(L(x)) \times \PP(L(x)^{\perp}) \right)$ lies in the
hyperplane $Z \subset \PP^{n-1}$ cut out by the equation $z_1+\cdots+z_n=0$. Let $\Gamma$ be the closure of the
graph of $m$ in $\PP(L(x)) \times \PP(L(x)^{\perp}) \times Z$, let
$\tilde{\Gamma}$ be a resolution of singularities of $\Gamma$ and let
$\tilde{m}: \tilde{\Gamma} \to Z$ be the composite map. Note that $\tilde{\Gamma}$ is $(d-1)+(n-d-1)=n-2$ dimensional, as is $Z$. 

Recall the definition of $\Omega_i$ for $i \leq n-1$: let $M$ be a generic $n-i$ plane
in $\CC^n$ and $\ell$ a generic line in $M$. Then $\Omega_i$ is the
set of $x \in G(d,n)$ such that $\ell \subset L(x)$ and $L(x) + M \neq
\CC^n$. In this section, we will write $\Omega_i(\ell, M)$ in order to record the dependence on $\ell$ and $M$.

\begin{prop}
  Suppose that $1 \leq i \leq n-1$. The variety $\left( \overline{T x}
    \cap \Omega_i(\ell,M) \right) \times \PP^{c-1} \times
  \PP^{\max(i-d,0)}$ is birational to $\tilde{m}^{-1}(W)$ where $W$ is a certain $\PP^{i-1}$ linearly embedded in $Z$. The linear space $W$ depends on
  $(\ell,M)$ and, if $(\ell,M)$ is chosen generically, then $W$ is generic
  in the Grassmannian of $\PP^{i-1}$'s in $Z$. Moreover, for $W$
  chosen generically, $\tilde{m}^{-1}(W)$ is smooth. 
\end{prop}

\begin{proof}
  We first prove this result in the case $c=1$ and $i \leq d$. The open subvariety $T x \cap
  \Omega_i(\ell,M)$ is dense in $ \overline{T x} \cap
  \Omega_i(\ell,M)$, so it is enough to understand $T x \cap
  \Omega_i(\ell,M)$. When $c=1$ we have $T x \isomorph T$ and we will describe
  $T x \cap \Omega_i(\ell,M)$ as a subvariety of $T$. Let $t$ be a point of $T$; the line $\ell$ is contained in $t \cdot L(x)$ if and only if $t^{-1} \cdot [\ell] \in
  \PP(L(x))$.  (Here $[\ell]$ is the class of $\ell$
  in $\PP^{n-1}$.) Note that $t^{-1} \cdot [\ell]$ is always in $(\CC^*)^{n-1} \subset \PP^{n-1}$.  Thus, we
  can think of $T x \cap \Omega_i(\ell,M)$ as a subvariety of
  $\PP(L(x)) \cap (\CC^*)^{n-1}$ by identifying $t \cdot L(x)$ with $t^{-1} [ \ell]$.

Our goal is to understand which
  points $u \in \PP(L(x)) \cap (\CC^*)^{n-1}$ correspond to $t$ such that
  $t \cdot L(x) + M \neq \CC^n$. We can restate this condition as $t^{-1} \cdot
  \PP(L(x)^{\perp}) \cap \PP(M^{\perp}) \neq \emptyset$. As $\ell$ is
  generic, all of its coordinates are nonzero and we may think of
  $\ell$ as a point of $T$. In this sense, the relation between $u$ and $t$ is $t=\ell \cdot u^{-1}$. Let $u \in \PP(L(x)) \cap (\CC^*)^{n-1}$,
  we want to understand when $(\ell^{-1} u) \cdot \PP(L(x)^{\perp}) \cap
  \PP(M^{\perp}) \neq \emptyset$. This happens if there is some $v \in
  \PP(L(x)^{\perp})$ such that $(\ell^{-1} u) \cdot v \in \PP(M^{\perp})$ or,
  equivalently, if $m(u,v) \in \ell \cdot \PP(M^{\perp})$.  When $i \leq d$ and $u$ corresponds to a generic point of $T x \cap \Omega_i(\ell,M)$, this point $v$ is unique.
  So $T x \cap \Omega_i(\ell,M)$ is birational to $m^{-1}(\ell \cdot
  \PP(M^{\perp}) \cap (\CC^*)^{n-1})$ and $ \overline{T x} \cap
  \Omega_i(\ell,M)$ is birational to $\tilde{m}^{-1}(\ell \cdot
  \PP(M^{\perp}))$. Clearly, $\ell\cdot \PP(M^{\perp})$ is an $(i-1)$
  dimensional projective space and, if $\ell$ and $M$ are chosen
  generically, $\ell \cdot \PP(M^{\perp})$ is a generic such space within $Z$.
  
  When $c>1$ the argument is basically the same except that $x$ has a
  nontrivial stabilizer in the $T$ action. Let $K \subset T$ be this
  stabilizer, we have $\dim K=c-1$. The torus $K$ acts on $\PP(L) \times
  \PP(L^{\perp})$ and a similar argument to the above shows that
  $m^{-1}(\ell \PP(M^{\perp}) \cap (\CC^*)^{n-1})/K$ is birational to $T x \cap
  \Omega_i(\ell,M)$. This exhibits $m^{-1}(\ell \cdot \PP(M^{\perp}) \cap
  (\CC^*)^{n-1})$ as a principal $K$-bundle over $T x \cap
  \Omega_i(\ell,M)$ and this bundle can be trivialized over some dense
  open subset $U$ of $T x \cap \Omega_i(\ell,M)$. So an open subset of
  $m^{-1}(\ell \cdot \PP(M^{\perp}))$ is isormorphic to $U \times K$ and thus
  to an open subset of $\overline{T x} \cap \Omega_i(\ell,M) \times
  \PP^{c-1}$.
  
  Similarly, when $i > d$, $v$ is no longer unique but, rather, the space of possible $v$'s is generically
  a $\PP^{i-d}$. Once again, it is easy to show that this bundle may
  be trivialized over an open set.

The smoothness of $\tilde{m}^{-1}(W)$ for generic $W$ follows from the Kleiman-Bertini theorem (see \cite{Kleim}) and the smoothness of $\tilde{\Gamma}$.
\end{proof}

\begin{cor} \label{BiRat}
  With the notations above, $\chi(\overline{T x} \cap \Omega_i) =
  \chi(\tilde{m}^{-1}(W))$ whenever $\Omega_i$ is chosen with respect
  to a generic $(\ell,M)$ and $W$ is a generic $(i-1)$-dimensional
  projective space in $Z$.
\end{cor}

\begin{proof}
  When $c=1$ and $i \leq d$, this is a direct consequence of the
  above. In general, this follows because $\chi(X \times Y)=\chi(X)
  \chi(Y)$ and $\chi(\PP^r)=1$.
\end{proof}

\textbf{Example:} We consider the examples of two points in $G(2,4)$. Our first example is a point whose matroid is the uniform matroid of rank $2$ on $4$ elements. We can take
$$L = \Span \begin{pmatrix} 1 & 0 & a & b \\ 0 & 1 & c & d \end{pmatrix} \quad L^{\perp} = \Span \begin{pmatrix} -a & -c & 1 & 0 \\ -b & -d & 0 & 1 \end{pmatrix}$$
where $abcd(ad-bc) \neq 0$. Then $\PP(L) \isomorph \PP(L^{\perp}) \isomorph \PP^1$, $Z \isomorph \PP^2$ and the map $m: \PP(L) \times \PP(L^{\perp}) \dashrightarrow Z$ can be factored as $\PP^1 \times \PP^1 \into \PP^3 \onto \PP^2$ where the first map is the Segre embedding and the second map is the linear projection away from $(1/a : -1/b : -1/c : 1/d)$. The assumption that $ad-bc \neq 0$ tells us that this point is not on $\PP^1 \times \PP^1$ so the composition is well defined and we can take $\tilde{\Gamma}=\PP^1 \times \PP^1$. The map $\tilde{\Gamma}=\PP^1 \times \PP^1 \to Z$ is of degree $2$, so the inverse image of a generic point of $Z$ is two points. The inverse image of a generic line in $Z$ is a $(1,1)$ curve in $\PP^1 \times \PP^1$ and hence has genus $0$ and holomorphic Euler characteristic $1$. The inverse image of $Z$ is $\PP^1 \times \PP^1$, which also has holomorphic Euler characteristic $1$. So, in this case, $r_{\overline{T x}}=2t+t^2+t^3+t^4 + \cdots = \frac{2t-t^2}{1-t}$ and $g_x(t)=2t+t^2$.

We now see what happens if we take $ad-bc=0$, but $abcd \neq 0$. This corresponds to a series-parallel matroid -- specifically, the parallel extension of the uniform matroid of rank $2$ on $3$ elements. Now the rational map $m : \PP^1 \times \PP^1 \dashrightarrow Z$ is given by projection from a point on $\PP^1 \times \PP^1$. Take $\tilde{\Gamma}$ to be the blow up of $\PP^1 \times \PP^1$ at that point. Now the map $\tilde{\Gamma} \to Z$ is only degree $1$. It is still true that the inverse image of a generic line in $Z$ has genus $0$ and the inverse image of $Z$ has holomorphic Euler characteristic $1$. So $r_{\overline{Tx}} = t+t^2+t^3+t^4 + \cdots = \frac{t}{1-t}$ and now $g_x(t)=t$.

At this point, we can prove three nontrivial results about $g_x(t)$.

\begin{prop} \label{Duality}
We have $g_x(t)=g_{x^{\perp}}(t)$.
\end{prop}

\begin{proof}
  $g_x(t)$ is expressed in terms of $\chi \left( \overline{T x} \cap
    \Omega_i \right)$ for various $i$ so it is enough to show that
  this quantity is invariant under exchanging $x$ and $x^{\perp}$. Our
  description of $\tilde{m}$ is symmetric under the exchange of $x$ and $x^{\perp}$, so this is clear.
\end{proof}

\begin{prop} \label{ValueAtMinusOne}
Writing $\dim \overline{Tx}=n-c$, we have $g_x(-1)=(-1)^c$.
\end{prop}

\begin{proof}
By definition, $(-1)^c g_x(-t)/(1-t) = \sum_i \chi(\overline{Tx} \cap \Omega_{i}) t^i$, so $(-1)^c g(-1)$ is the residue of $\sum_i \chi(\overline{Tx} \cap \Omega_{i}) t^i$ at $t=1$. We will show that, for $i$ sufficiently large, $ \chi(\overline{Tx} \cap \Omega_{i}) =1$. For $i \geq n-1$, $\overline{Tx} \cap \Omega_{i}$ is, by definition, $\overline{Tx} \cap \Omega_{n-1}$. We have just shown that $\chi(\overline{Tx} \cap \Omega_{n-1})=\chi(\tilde{m}^{-1}(\PP^{n-2}))=\tilde{\Gamma}$. But $\tilde{\Gamma}$ is birational to $\PP(L(x)) \times \PP(L(x)^{\perp})$ and thus has holomorphic Euler characteristic $1$.
\end{proof}

For the next result, we need to introduce some notation.  Let $e \in
[n]$. Let $x \in G(d,n)$ and assume that $x$ is not in any coordinate
subgrassmannian. Let $p_e$ be the map $\CC^n \into \CC^{n+1}$ by
$(u_1,\ldots, u_n) \mapsto (u_1, \ldots, u_e, \ldots, u_n, u_e)$ and
let $p_e(x)=L^{-1}(p_e(L(x))$. Let $s_e(x)=(p_e(x^{\perp}))^{\perp}$. We call
$p_e(x)$ and $s_e(x)$ the parallel extension and series coextension (respectively) of
$x$ at $e$. We could also define $s_e(x)$ by $L(s_e(x))=\iota(L(x)) \oplus
\CC(e_{n+1}-e_e)$ where $\iota$ is the embedding of $\CC^n$ into the
first $n$ coordinates of $\CC^{n+1}$. The matroids $M(p_e(x))$ and
$M(s_e(x))$ depend only on $M(x)$ and $e$. These matroids are called
the parallel extension and series coextension of $M(x)$ at $e$. See section 7.6 of \cite{White1} for more background on these constructions, which are described there as special cases of the more general operations of series and parallel connection.

\begin{prop} \label{SerInvariant}
We have $g_x=g_{s_e(x)}=g_{p_e(x)}$. 
\end{prop}

This proposition explains why it a good idea to have the sum
defining $g$ stretch out to infinity, it would not be true if we
truncated the sum at $d$, $n$ or some other natural point.

\begin{proof}
 We first show that $g_{p_e(x)}=g_x$.  We will write $\Omega'_i$, $T'$, \emph{etc.} to denote objects associated with $p_e(x)$.
  We must show that, for every positive integer $i$,
  $\chi(\overline{T' p_e(x)} \cap \Omega'_i) = \chi(\overline{T x} \cap
  \Omega_i)$. We first show this in the case $i \leq n-1$.
  
  For this purpose, we do not use corollary \ref{BiRat} but work directly
  with the variety $\overline{T' p_e(x)} \cap \Omega'_i(\ell',M')$. We
  note that, for $\ell'$ and $M'$ generic, the open subvariety $T' p_e(x) \cap
  \Omega'_i(\ell',M')$ is dense in $\overline{T' p_e(x)} \cap
  \Omega'_i(\ell',M')$. If $\ell' \in t' L(p_e(x))$ then we must have
  $t'_{n+1}/t'_e=\ell'_{n+1}/\ell'_e$ (where $t'=(t_1, \ldots,
  t_{n+1})$ and $\ell'=(\ell_1,\ldots, \ell_{n+1})$.) The set of $t' \in T'$ with
  this property is a principle homogeneous space for $T$ and,
  identifying it with $T$, we get an isomorphism between $T' p_e(x)
  \cap \Omega'_i(\ell',M')$ and $T x \cap \Omega_i(\ell,M)$. Here
  $\ell$ is the projection of $\ell'$ onto the first $n$ coordinates
  and $M$ is the projection onto the first $n$ coordinates of $M' \cap
  \{ \ell'_e x_{n+1}-\ell'_{n+1} x'_e=0 \}$. Since $T' p_e(x) \cap
  \Omega'_i(\ell',M')$ and $T x \cap \Omega_i(\ell,M)$ are isomorphic,
  $\overline{T' p_e(x)} \cap \Omega'_i(\ell',M')$ and $\overline{T x}
  \cap \Omega_i(\ell,M)$ are birational. It is easy to see that $\ell
  \subset M$ and, if $(\ell',M')$ was generic, so is $(\ell,M)$.
  
  We now consider the case of $i \geq n$. We have
  $\Omega_{n-1}=\Omega_n=\Omega_{n+1}=\cdots$ and
  $\Omega'_{n}=\Omega'_{n+1}=\cdots$ by definition. Thus, we will be
  done if we can show that $\chi(\overline{T p_e(x)} \cap
  \Omega'_n)=\chi(\overline{T x} \cap \Omega_{n-1})$. These two
  quantities are both $1$ as we observed in the proof of the previous proposition. We have now shown that $g_{p_e(x)}=g_x$.

Using proposition \ref{Duality} twice, we have $g_{s_e(x)}=g_{s_e(x)^{\perp}}=g_{p_e(x^{\perp})}=g_{x^{\perp}}=g_x$. So $g_{s_e(x)}=g_x$.
\end{proof}

A matroid is called series-parallel if it can be obtained by repeated series-parallel extensions from the matroid corresponding to a generic point in $G(1,2)$. See section 6.4 of \cite{White1} for background on series-parallel matroids. The following corollary logically belongs in the next section, but it fits more naturally here -- the reader can check that no circularity is involved.

\begin{cor} \label{SerParGivesX}
Let $x \in G(d,n)$ and assume that $x$ is not in any coordinate subgrassmannian. Then $M(x)$ is series-parallel if and only if $g_x(t)=t$. $M(x)$ is a direct sum of $c$ series-parallel matroids if and only if $g_x(t)=t^c$.
\end{cor}

\begin{proof}
First, assume that $x$ is series-parallel; we want to show that $g_x(t)=t$. By proposition \ref{SerInvariant},  it is enough to consider the case that $(d,n)=(1,2)$ and $x$ is a generic point of $G(1,2)=\PP^1$. In this case, $\PP(L(x))$, $\PP(L(x)^{\perp})$ and $Z$ are all points so $\tilde{m}^{-1}(\PP^{\min(i,0)})$ is a point for all $i$ and $\chi(\tilde{m}^{-1}(\PP^{\min(i,0)}))=1$. Then $h_x(t) = (1-t) \sum_{i=1}^{\infty} t^i=t$ and $g_x(t)=-h_x(-t)=t$. The case where $M(x)$ is a direct sum of $c$ series-parallel matroids then follows easily from proposition \ref{ProductFormula}.

For the converse, suppose that $g_x(t)=t$. Then, by theorem \ref{Beta}, $\beta(M(x))=1$. By theorem 7.6 of \cite{Bry}, this implies that $M(x)$ is series-parallel. Similarly, suppose that $g_x=t^c$. Suppose that $M(x)$ has $r$ connected components, so $x=x_1 \oplus \cdots \oplus x_r$ with $x_i$ connected. Because $x$ is not contained in any coordinate subgrassmannian, each connected component of $M(x)$ has at least $2$ elements.The coefficient of $t$ in $g_{x_i}(t)$ is nonzero and the constant term of $g_{x_i}(t)$ is zero. So $r$ is precisely the power of $t$ that divides $g_x(t)=\prod g_{x_i}(t)$ and we have $r=c$. Moreover, we must have $g_{x_i}(t)=t$ for each $i$. Then, as before, each $M(x_i)$ is a series-parallel matroid. 
\end{proof}

\section{Proof of proposition \ref{ProductFormula}} \label{ProdProof}

In this section, we will prove proposition \ref{ProductFormula}, which states:

\begin{PF}
If $[n]=\bigsqcup_{k=1}^r A_k$ and $x_k$ is a point of $G(d_k,A_k)$ not contained in any coordinate subgrassmannian then $g_{\bigoplus_{k=1}^r x_k}=\prod_{k=1}^r g_{x_k}$.
\end{PF}

\begin{proof}
Clearly, it is enough to prove the result in the case $r=2$. Let $n_k=|A_k|$ for $k=1$, $2$ and let $L_k:=L(x_k)$, $L_k^{\perp}:=L(x_k)^{\perp}$, $m_k$, $\tilde{\Gamma}_k$, $\tilde{m}_k$, $Z_k$ and so forth have the obvious meanings with respect to $x_k$. Let $L$, $L^{\perp}$, $m$, $\tilde{m}$, $\tilde{\Gamma}$, $Z$ and so forth have the corresponding meanings with respect to $x$. Let $n_k - c_k = \dim( (\CC^*)^{A_k} x_k)$ and $n-c=\dim ((\CC^*)^n x)$, so $c=c_1+c_2$.

Our goal is to establish the equality
\begin{multline*}
(-1)^{c_1+c_2-1} (1-t) \left( \sum_i \chi( \tilde{m}^{-1}(\PP^{\min(i-1,n-2)})) t^i \right) = \\
(-1)^{c_1-1} (1-t)  \left( \sum_{i_1} \chi( \tilde{m_1}^{-1}(\PP^{\min(i_1-1,n_1-2)})) t^{i_1} \right) \times \\
(-1)^{c_2-1} (1-t)  \left( \sum_{i_2} \chi( \tilde{m_2}^{-1}(\PP^{\min(i_2-1,n_2-2)})) t^{i_2} \right)
\end{multline*}
or, equating coefficients of $t^i$,
\begin{multline}
\chi(\tilde{m}^{-1}(\PP^{\min(i-1,n-2)})) = \\
\sum_{i_1+i_2=i}  \chi( \tilde{m_1}^{-1}(\PP^{\min(i_1-1,n_1-2)}))  \chi( \tilde{m_2}^{-1}(\PP^{\min(i_2-1,n_2-2)})) - \\
\sum_{i_1+i_2=i-1}  \chi( \tilde{m_1}^{-1}(\PP^{\min(i_1-1,n_1-2)}))  \chi( \tilde{m_2}^{-1}(\PP^{\min(i_2-1,n_2-2)})) \label{showthis0}
\end{multline}
In these equations, $\PP^j$ should always be interpreted as a generic $\PP^j$ in $Z$, $Z_1$ or $Z_2$ as appropriate.

Fix a value of $i$ for which we will establish (\ref{showthis0}). Both sides of (\ref{showthis0}) are invariant under series and parallel extension of $x_1$ and $x_2$; by making enough such extensions we can assume that $d_1$, $d_2$, $n_1-d_1$ and $n_2-d_2$ are all greater than $i+1$. As a consequence, all of the $\min$'s in equation (\ref{showthis0}) drop out. 

Let $W \subset Z$ be the hyperplane where $\sum_{a \in A_1} z_a=0$
(equivalently $\sum_{a \in A_2} z_a=0$). Then there is a rational map
$q : W \dashrightarrow Z_1 \times Z_2$ -- specifically, $q$ is the quotient of the
obvious isomorphism $\CC^{A_1 \sqcup A_2} \to \CC^{A_1} \times
\CC^{A_2}$ by the actions of $\CC^*$ acting on $\CC^{A_1 \sqcup A_2}$ and $(\CC^*)^2$ acting on $\CC^{A_1} \times \CC^{A_2}$. Let $U \subset W$ be the open locus on which $q$ is
defined, $W \setminus U$ has codimension $\min(n_1-1,n_2-1)$. The map $q:U \to
Z_1 \times Z_2$ is a $\CC^*$ bundle. Similarly, there are rational
maps $r: \PP(L) \to \PP(L_1) \times \PP(L_2)$ and $r^{\perp} :
\PP(L^{\perp}) \to \PP(L_1^{\perp}) \times \PP(L_2^{\perp})$; let $V$
and $V^{\perp}$ be the loci where $r$ and $r^{\perp}$ are
defined. Then $\PP(L) \setminus V$ has codimension $\min(d_1-1,d_2-1)$ and
$\PP(L^{\perp}) \setminus V^{\perp}$ has codimension
$\min(n_1-d_1,n_2-d_2)$.  

As rational maps, we have $q \circ m = (m_1
\times m_2) \circ (r \times r^{\perp})$. We can extend $r \times
r^{\perp}$ to a rational map $\tilde{s} : \tilde{\Gamma} \dashrightarrow \tilde{\Gamma_1}
\times \tilde{\Gamma_2}$. By altering out choice of $\tilde{\Gamma}$, we may assume that $\tilde{s}$ is a well defined morphism. Over a generic point of $\tilde{\Gamma}_1 \times \tilde{\Gamma_2}$, the fiber of $\tilde{s}$ is some compactification of $(\CC^*)^2$ (which one depends on the choice of resolution of singularities $\tilde{\Gamma}$). Then, in a slight abuse of notation, $q \circ \tilde{m} =
(\tilde{m_1} \times \tilde{m_2}) \circ \tilde{s}$ where the right hand
side is a well defined morphism but the left hand side is only a
rational map.

Our goal is to compute the holomorphic Euler characteristic of
$\tilde{m}^{-1}(\PP^{i-1})$ where $\PP^{i-1}$ is chosen generically in
$Z$. The image of $\tilde{m}$ lies in $W$, so we must compute
$\tilde{m}^{-1}(\PP^{i-1} \cap W)$. We will denote $\PP^{i-1} \cap W$
by $H$; $H$ is a generic $(i-2)$-plane in $W$. As we took $n_1$ and
$n_2$ large, we may assume that $q(H)$ is well defined and isomorphic
to $H$, that the projections of $q(H)$ to $Z_1$ and $Z_2$ are
$\PP^{i-2}$'s linearly embedded in $Z_1$ and $Z_2$ and that $q(H)$ is the
graph of an isomorphism between these $\PP^{i-2}$'s. Let us call a
subvariety $K$ of $Z_1 \times Z_2$ a diagonal $\PP^{i-2}$ if $K$ is
the graph of an isomorphism between a $\PP^{i-2}$ linearly embedded in
$Z_1$ and a $\PP^{i-2}$ linearly embedded in $Z_2$. The group $\PGL(Z_1) \times
\PGL(Z_2)$ acts transitively on the collection of diagonal
$\PP^{i-2}$'s and the reader may easily check that, if $\PP^{i-1}$ is
chosen generically in $Z$, then $q(W \cap \PP^{i-1})$ is a generic
diagonal $\PP^{i-2}$.

We claim that $(\tilde{s}^{-1} \circ (\tilde{m_1} \times \tilde{m_2})^{-1})(
q(H)) =\tilde{m}^{-1}(\PP^{i-1})$. Note that the left hand side
contains the right as the image of $\tilde{m}$ lands in
$W$. Moreover, the inverse image of the open locus in $H$ where all
$n$ coordinate functions are nonzero is dense in both the left and
right hand side. Therefore, to see that the two sides are equal, it is
enough to see that both are smooth varieties. On the right hand side,
we know that $\tilde{\Gamma}$ is smooth, $\PP^{i-1}$ is smooth and
$\PP^{i-1}$ was chosen generically under the $\PGL(Z)$ action on $Z$,
so $\tilde{m}^{-1}(\PP^{i-1})$ is smooth by the Kleiman-Bertini theorem (see \cite{Kleim}). Similar
results apply to the left hand side, using the smoothness of $q(H)$
and the $\PGL(Z_1) \times \PGL(Z_2)$ action.  So we may concentrate on
computing $\chi(\tilde{s}^{-1} (\tilde{m_1} \times \tilde{m_2})^{-1}
q(H) )$. Moreover there is a dense open subset of
$(\tilde{s}^{-1} \circ (\tilde{m_1} \times \tilde{m_2})^{-1})( q(H))$ which is a
$(\CC^*)^2$ bundle over a dense open subset of $(\tilde{m_1} \times
\tilde{m_2})^{-1}( q(H))$ so we conclude that $(\tilde{s}^{-1} \circ
(\tilde{m_1} \times \tilde{m_2})^{-1})( q(H))$ is birational to $\PP^2
\times (\tilde{m_1} \times \tilde{m_2})^{-1} q(H)$. The latter variety
is clearly proper and we may use the Kleiman-Bertini theorem applied to
$\tilde{m}_1 \times \tilde{m}_2$ to
conclude that it is smooth. Thus, $\chi((\tilde{s}^{-1} \circ (\tilde{m_1}
\times \tilde{m_2})^{-1})( q(H)) )=\chi((\tilde{m_1} \times
\tilde{m_2})^{-1} (q(H)))$.

Our goal now is to show that
\begin{multline}
\chi((\tilde{m_1} \times \tilde{m_2})^{-1}(q(H))) = \\
\sum_{i_1+i_2=i}  \chi( \tilde{m_1}^{-1}(\PP^{i_1-1}))  \chi( \tilde{m_2}^{-1}(\PP^{i_2-1})) - 
\sum_{i_1+i_2=i-1}  \chi( \tilde{m_1}^{-1}(\PP^{i_1-1}))  \chi( \tilde{m_2}^{-1}(\PP^{i_2-1})). \label{showthis1}
\end{multline}

In lemma \ref{KEquality}, we show that, in $K^{\circ}(Z_1 \times Z_2)$, we have
\begin{equation}
[\OO_{q(H)}] = \sum_{i_1+i_2=i}  [\OO_{\PP^{i_1-1} \times \PP^{i_2-1}}] - 
\sum_{i_1+i_2=i-1}  [\OO_{\PP^{i_1-1} \times \PP^{i_2-1}}] . \label{showthis2}
\end{equation}
Assuming this, we may pull this equality back along $\tilde{m}_1 \times \tilde{m_2}$ to get an equality in $K^{\circ}(\tilde{\Gamma_1} \times \tilde{\Gamma}_2)$. In general, the formula for pullback involves higher $\TorTor$'s, but we may use lemma \ref{BertiniTor} and the transitive action of $\PGL(Z_1) \times \PGL(Z_2)$ to assume that all the higher $\TorTor$'s drop out. Then applying $\chi$ to the equality in $K^{\circ}(\tilde{\Gamma_1} \times \tilde{\Gamma}_2)$ (and using $\chi(A \times B)=\chi(A) \chi(B)$) yields equation (\ref{showthis2}) and we are done.
\end{proof}

\begin{lemma} \label{KEquality}
Let $q(H)$ be a diagonal $\PP^{i-2}$ in $Z_1 \times Z_2$. Then, in $K^{\circ}(Z_1 \times Z_2)$, we have
$$[\OO_{q(H)}] = \sum_{i_1+i_2=i}  [\OO_{\PP^{i_1-1} \times \PP^{i_2-1}}] - 
\sum_{i_1+i_2=i-1}  [\OO_{\PP^{i_1-1} \times \PP^{i_2-1}}] .$$
\end{lemma}

\begin{proof}
Let $\pi_1$ and $\pi_2$ be the projections of $Z_1 \times Z_2$ onto its factors. We may assume that $\dim Z_1 =\dim Z_2=i-2$, as otherwise we can first prove the equality in $\pi_1(q(H)) \times \pi_2(q(H))$ and then push it forward along the closed inclusion $\pi_1(q(H)) \times \pi_2(q(H)) \into Z_1 \times Z_2$. Also, by changing coordinates on $Z_1$ and $Z_2$, we may assume that $q(H)$ is the diagonal in $Z_1 \times Z_2 \isomorph \PP^{i-2} \times \PP^{i-2}$. 

Thus, our goal is to compute the class of the diagonal in $K^{\circ}(\PP^{i-2} \times \PP^{i-2})$. This can be done in many ways, we simply cite \cite{Brion} which gives a formula for the class of the diagonal in $G/P$ for $G$ any reductive Lie group and $P$ a parabolic subgroup.
\end{proof}

\section{Proof of proposition \ref{Positivity}} \label{PosProof}

In this section, we will prove proposition \ref{Positivity}. This
result is indespensible in proving the Main Theorem, and is the part
of the paper which uses the most deep algebraic geometry. In
particular, while I am quite confident that the characteristic zero
hypotheses, and perhaps even the realizability hypotheses, could be
removed from the rest of this paper, I am uncertain as to whether this
result will remain true in positive characteristic.

We recall the result we are trying to prove, using corollary \ref{BiRat} to rewrite our statement:

\begin{PT}
  Let $x \in G(d,n)$ and assume $x$ is not contained in any coordinate
  subgrassmannian. Let $n-c = \dim Tx$. Let $\tilde{m} :
  \tilde{\Gamma} \to Z$ be as discussed above. Define $g_x(t)$ by 
$$\frac{(-1)^c g_x(-t)}{1-t} = \sum_{j=1}^{\infty} \chi(\tilde{m}^{-1}(\PP^{\min(j-1,n-2)})) t^j$$
where $\PP^{j-1}$ denotes a generic $\PP^{j-1}$ in $Z$. Then $g_x(t)$ has nonnegative coefficients, and the coefficient of $t^c$ is positive.
\end{PT}

\begin{proof}
  First, by proposition \ref{ProductFormula}, we can reduce to the
  case $c=1$. Also, in this case, we showed in theorem \ref{Beta} that
  the coefficient of $t$ is nonnegative in this circumstance. If $j
  \geq n$, then it is clear that the coefficient of $t^j$ is zero.
  (Actually, the coefficient of $t^j$ is zero if $j > \min(d,n-d)$.) So we will concentrate on showing the
  coefficient of $t^j$ is nonnegative when $2 \leq j < n$. In this case, we are being asked to show that
$$(-1)^{j-1} \left( \chi(\tilde{m}^{-1}(\PP^{j-1})) - \chi(\tilde{m}^{-1}(\PP^{j-2})) \right) \geq 0.$$
Now, on $Z$, we have the short exact sequence of sheaves
$$0 \to \OO_{\PP^{j-1}}(-1) \to \OO_{\PP^{j-1}} \to \OO_{\PP^{j-2}}
\to 0.$$
By proposition \ref{BertiniTor}, $\TorTor^{\tilde{m}}_1(
\OO_{\PP^{j-2}}, \OO_{\tilde{\Gamma}})=0$ for a generic choice of $\PP^{j-2}$ so this
sequence remains exact after pullback to $\tilde{\Gamma}$.

So we are being asked to show that
$$0 \leq (-1)^{j-1} \chi(\tilde{m}^* \OO_{\PP^{j-1}}(-1))$$
or, explicitly,
$$0 \leq (-1)^{j-1} \sum_{i=0}^{j-1} (-1)^i \dim H^i(\tilde{m}^*
(\OO_{\PP^{j-1}}(-1)), \tilde{m}^{-1}(\PP^{j-1})) .$$

We claim that  $\tilde{m} :
\tilde{\Gamma} \to Z$ is surjective and generically finite. Proof:
$\tilde{\Gamma}$ and $Z$ have the same dimension ($n-2$) so it is
enough to show that the fiber over a generic point of $Z$ is nonempty. By proposition
\ref{Birational}, this fiber is birational  to $\overline{Tx}
\cap \Omega_1(\ell,H)$ for a generic pair $(\ell,H)$ of ``line contained in hyperplane".  By theorem
\ref{Beta}, $\overline{Tx} \cap \Omega_1(\ell,H)$ consists of
$\beta(M(x))$ points and, as we noted in proposition \ref{CoeffOfT},  $\beta(M(x))>0$ when $c=1$.

Since we have chosen $\PP^{j-1}$
generically, $\tilde{m}^{-1}(\PP^{j-1})$ is smooth and $\tilde{m}$
restricted to $\tilde{m}^{-1}(\PP^{j-1})$ is surjective and
generically finite as a map to $\PP^{j-1}$. We abbreviate the line bundle
$\tilde{m}^*(\OO(-1))|_{\tilde{m}^{-1}(\PP^{j-1})}$ by $L$. Then $L$
has positive degree and $L$ restricted to any curve in
$\tilde{m}^{-1}(\PP^{j-1})$ has nonnegative degree, \emph{i.e.}, $L$
is nef. So, by Kawamata-Viehweg vanishing (\cite{Kawa}, \cite{Vieh}), $H^i(L,
\tilde{m}^{-1}(\PP^{j-1}))=0$ for $i \neq j-1$. Thus, the quantity we
are being asked to show is nonnegative is 
$$(-1)^{2j-2} \dim H^{j-1}(L, \tilde{m}^{-1}(\PP^{j-1}))=\dim
H^{2j-2}(L, \tilde{m}^{-1}(\PP^{j-1})).$$
Of course, the dimension of
a vector space is nonnegative, so we are done.
\end{proof}

\section{Behavior under $2$-sum} \label{SectionTwoSum}

In section, we will prove a result describing the behavior of $g_x$ under an operation called ``$2$-sum".  While this result is not used in proving our Main Theorem, it is invaluable in computing $g_x$ in practice. Let $A_1$ and $A_2$ be finite sets, let $x_r \in G(d_r,A_r)$ for $r=1$, $2$ and let $e_r \in A_r$. Assume that neither of the  $x_r$ is contained in any coordinate subgrassmannian. Let $\alpha \in \CC^*$. We define a point $x_1 +^{e_1,e_2}_2(\alpha) x_2 \in G(d_1+d_2-1,A_1 \sqcup A_2 \setminus \{ e_1,e_2 \})$ as follows: $L(x_1 +_2^{e_1,e_2}(\alpha) x_2)$ is the projection onto $\CC^{A_1 \sqcup A_2 \setminus \{ e_1,e_2 \}}$ of $\left( L(x_1) \oplus L(x_2) \right) \cap \{ z_{e_1} = \alpha z_{e_2} \}$. The points $x_1 +_2^{e_1,e_2}(\alpha) x_2$ and $x_1 +_2^{e_1,e_2}(\alpha') x_2$ lie in the same $T$-orbit for any two values $\alpha$ and $\alpha' \in \CC^*$, so we see that the matroid of $x_1 +^{e_1,e_2}_2(\alpha) x_2$ is independent of $\alpha$. We will therefore drop the $\alpha$ from our notation when dealing with quantities that only involve the $T$-orbit closure or only involve the matroid. The matroid $M(x_1 +_2^{e_1,e_2} x_2)$ is traditionally denoted $M(x_1) +_2 M(x_2)$ and called the $2$-sum of $M(x_1)$ and $M(x_2)$. This is an abuse of notation, as the matroid depends not only on $M(x_1)$ and $M(x_2)$ but also on $e_1$ and $e_2$. We will usually follow the matroid convention and drop the superscripted $e_1$ and $e_2$ from our notation. The reader should observe that $(x_1 +_2(\alpha) x_2)^{\perp} = x_1^{\perp} +_2(-\alpha^{ -1}) x_2^{\perp}$. See \cite{CE} and section 7.6 of \cite{White1} for more on this operation.

We will spend the rest of this section proving the following result:
\begin{prop}
With the above notation, we have
$$g_{x_1 +_2 \ x_2}(t) = g_{x_1}(t) g_{x_2}(t)/t.$$
\end{prop}

Before beginning our proof, we
fix some notations: Let $n_r$ be the cardinality of $A_r$ and let
$n_r-c_r$ be the dimension of $(\CC^*)^{A_r} x_r$. Let
$n=n_1+n_2-2$ and $n'=n_1+n_2$, $c=c_1+c_2-1$ and $c'=c_1+c_2$,
$d=d_1+d_2-1$ and $d'=d_1+d_2$. Let $A=A_1 \sqcup A_2 \setminus \{
e_1, e_2 \}$ and $A'=A_1 \sqcup A_2$. Let $x=x_1
+_2^{e_1,e_2}(\alpha) x_2$ and $x'=x_1 \oplus x_2$,
$L=L(x)$ and $L'=L(x')$. Let $Z_1$, $Z_2$, $Z$ and
$Z'$ be the hyperplanes in $\PP(\CC^{A_1})$, $\PP(\CC^{A_2})$,
$\PP(\CC^A)$ and $\PP(\CC^{A'})$ where the sum of the coordinates is
zero. In general, we use subscripts $1$ and $2$ to denote objects
associated $x_1$ and $x_2$, a lack of demarcation to denote objects associated
with $x=x_1 +^{e_1,e_2}_2(\alpha) x_2$ and primes to denote objects associated
with $x_1 \oplus x_2$. The meaning of symbols such as
$\tilde{\Gamma}$, $\tilde{m}$, \emph{etc.} should be
clear. We use $\pi_1$ for the projection $\tilde{\Gamma_1} \to \PP(L_1) \times \PP(L_1^{\perp})$ and define $\pi_2$, $\pi$ and $\pi'$ analogously.  Let $W' \subset Z'$ be the hyperplane where $\sum_{a \in A_1} z_a=0$ (equivalently, where $\sum_{a \in A_2} z_a=0$.)

\begin{proof}
First, it is easy to check that we have $(x \oplus y) +_2^{e_1,e_2}(\alpha) z = x \oplus (y  +_2^{e_1,e_2}(\alpha) z)$. Using this equality and proposition \ref{ProductFormula}, we immediately reduce to the case $c_1=c_2=1$. 

  By theorem \ref{ProductFormula}, our goal is to prove that $t g_{x_1
    +_2 \ x_2}(t) = g_{x_1 \oplus x_2}(t)$. This is
  equivalent to showing, for every $i \geq 1$, that $\chi(\tilde{m}^{-1}(\PP^{i-1}))= \chi((\tilde{m}')^{-1}(\PP^i))$ where $\PP^{i-1}$ and $\PP^{i}$ are
  generically chosen in $Z$ and $Z'$ respectively. Equivalently, we may show that $\chi(\tilde{m}^{\-1}(\PP^{i-1}))=\chi((\tilde{m}')^{-1}(\PP^{i-1}))$ where the first $\PP^{i-1}$ is chosen generically in $Z$ and the second is chosen generically in $W'$. We can identify $Z$ with the hyperplane $\{ z_{e_1} = - z_{e_2} \}$ in $W'$ in an obvious way; we write $\iota$ for the resulting injection $Z \into W'$.  Our proof breaks into two parts; first we show that $\chi(\tilde{m}^{-1}(\PP^{i-1}))=\chi((\tilde{m}')^{-1}(\PP^{i-1}))$ where the first $\PP^{i-1}$ is chosen generically in $Z$ and the second in $\iota(Z)$; then we will show that taking $\PP^{i-1}$ generic in $\iota(Z)$ instead of in all of $W$ does not change the generic value of $\chi(\tilde{m}^{-1}(\PP^{i-1}))$.

\textbf{Part 1:} Let $H$ be a generic $\PP^{i-1}$ in $Z$. We will show that $\tilde{m}^{-1}(H) \times \PP^1$ is birational to $(\tilde{m}')^{-1}(\iota(H))$ and that (after making an appropriate choice of $\tilde{\Gamma}'$) both are smooth and proper. This implies the equality of Euler characteristics which is our first goal. 

Clearly, if $H$ is chosen generically in $Z$ then $\iota(H)$ is chosen generically in $\iota(Z)$. We know that $\tilde{m}^{-1}(H)$ is smooth by the Kleiman-Bertini theorem (see \cite{Kleim}) as $\tilde{\Gamma}$ is smooth and $H$ is chosen generically. If we knew that $(\tilde{m}')^{-1}(\iota(Z))$ was smooth then the same argument would show that $(\tilde{m}')^{-1}(\iota(H))$ is smooth. Let $F \subset \PP^1 \times W'$ be the pencil of hyperplanes with $F_{a_1 : a_2} := \{ a_1 z_{e_1} = a_2 z_{e_2} \}$ over $(a_1 : a_2) \in \PP^1$; the fiber $F_{1:-1}$ is $\iota(Z)$. We may assume that $\tilde{m}' : \tilde{\Gamma}' \to W'$ factors through $F$ by the standard trick -- take the connected component of $\tilde{\Gamma}' \times_{W'} F$ which lies over the generic point of $\tilde{\Gamma}'$, resolve its singularities and use this resolution to replace $\tilde{\Gamma}'$. Once we have done this, by Kleiman-Bertini applied to $\tilde{\Gamma} \to \PP^1$, we know that $(\tilde{m}')^{-1}(F_{a_1:a_2})$ is smooth for a generic $(a_1:a_2) \in \PP^1$. Now, there is an action $\rho$ of $\CC^*$ on $\PP(L_1 \oplus L_2) \times \PP(L_1^{\perp} \oplus L_2^{\perp})$ by scaling the coordinates of $L_2$ and leaving alone those of $L_1$, $L_1^{\perp}$ and $L_2^{\perp}$; there is a similar action $\sigma$ of $\CC^*$ on $Z'$ which scales the coordinates indexed by $A_2$ and leaves alone those indexed by $A_1$. The rational map $m$ relates these actions in the sense that $m(\rho(t)(y))=\sigma(t)m(y)$ when both maps are defined. For any $t \in \CC^*$, we can define a new $\tilde{\Gamma}'$, which we will denote $\tilde{\Gamma}'_t$, by taking $\tilde{\Gamma}'_t$ abstractly isomorphic to $\tilde{\Gamma}'$ but replacing $\pi$ with $\pi_t := \rho(t^{-1}) \circ \pi$ and replacing $\tilde{m}'$ with $\tilde{m}'_t :=\sigma(t) \circ \tilde{m}'$.  Then $(\tilde{m}'_t)^{-1}(F_{a_1:t a_2}) \isomorph \tilde{m}^{-1}(F_{a_1:a_2})$. So, by replacing $\tilde{\Gamma}'$ by $\tilde{\Gamma}'_t$ for an appropriate $t$, we may assume that $(\tilde{m}')^{-1}(F_{1:-1}))=(\tilde{m}')^{-1}(\iota(Z))$ is smooth.

We now explain why $\tilde{m}^{-1}(H) \times \PP^1$ is birational to $(\tilde{m}')^{-1}(\iota(H))$. Let $U'$ be the open subset of $\PP(L') \times \PP((L')^{\perp})$ where the $e_1$ and $e_2$ coordinates are nonzero in each factor. Let $K$ be the hyperplane in $\PP(L)$ defined by the equation $\sum_{a \in A_1 \setminus \{ e_1 \}} x_a=0$ and define $K^{\perp}$ similarly. Let $U \subset \PP(L) \times \PP(L^{\perp})$ be the complement of $(K \times \PP(L^{\perp})) \cup (\PP(L) \times K^{\perp})$. Then
$\pi^{-1}(U) \cap \tilde{m}^{-1}(H)$ and $(\pi')^{-1}(U') \cap (\tilde{m}')^{-1}(\iota(H))$ are dense in $\tilde{m}^{-1}(H)$ and $(\tilde{m}')^{-1}(\iota(H))$ respectively, thus it suffices to show that $(\pi^{-1}(U) \cap \tilde{m}^{-1}(H)) \times \CC^*$ is birational to $(\pi')^{-1}(U') \cap (\tilde{m}')^{-1}(\iota(H))$. Now, the rational maps $m$ and $m'$ are well defined on $U$ and $U'$, so we may assume that $\pi^{-1}(U) \isomorph U$ and $(\pi')^{-1}(U') \isomorph U'$. We write $\mu$ and $\mu'$ for the restrictions of $m$ and $m'$ to $U$ and $U'$. We now see that showing $(\pi^{-1}(U) \cap \tilde{m}^{-1}(H)) \times \CC^*$ is birational to $(\pi')^{-1}(U') \cap (\tilde{m}')^{-1}(\iota(H))$ is equivalent to showing that $\mu^{-1}(H) \times \CC^*$ is birational to $(\mu')^{-1}(\iota(H))$. In fact, we will show that, under the reductions we have already made, $\mu^{-1}(H) \times \CC^*$  and  $(\mu')^{-1}(\iota(H))$ are isomorphic.

There is an action $\rho$ of $\CC^*$ on $U'$ where $\rho(t)$ scales $L_2$ by $t$, $L_2^{\perp}$ by $t^{-1}$ and leaves $L_1$ and $L_1^{\perp}$ alone. This action is free, the quotient is clearly identified with $U$ and each orbit contains exactly one point of $\PP(L' \cap \{x_{e_1} = \alpha x_{e_2} \}) \times  \PP((L')^{\perp} \cap \{x_{e_1} = -1/\alpha x_{e_2} \})=\PP(L) \times \PP(L^{\perp})$. We thus obtain a natural isomorphism $U' \isomorph U \times \CC^*$. Writing $p$ for the projection $U' \to U$, we have $\iota \circ \mu \circ p = \mu'$. Thus, $(\mu')^{-1}(\iota(H))$ is a trivial $\CC^*$ bundle over $\mu^{-1}(H)$, as promised.

\textbf{Part 2:} We now show that $\chi((\tilde{m}')^{-1}(H))$ is the same for $H$ a generic $\PP^{i-1}$ in $\iota(Z)$ or for $H$ a generic $\PP^{i-1}$ in $W'$. The point is the following: by the well definedness of pullback and $\chi$ as maps on $K^{\circ}$ we know that $\sum (-1)^i \chi(\TorTor_i([\OO_H], \tilde{\Gamma}'))$ is completely independent of $H$. For $H$ chosen generically in $W'$, we know that all the higher $\TorTor$'s vanish so it is enough to show that, for $H$ chosen generically in $\iota(Z)$, we still have this $\TorTor$ vanishing. This may be checked in two steps: first, we check that there are no higher $\TorTor$'s when pulling $[\OO_{\iota(Z)}]$ back to $\tilde{\Gamma}'$ and, second, we check that, for $H$ chosen generically in $\iota(Z)$, there are also no higher $\TorTor$'s when pulling $[\OO_H]$ back to $(\tilde{m}')^{-1}(\iota(Z))$. To see the first $\TorTor$ vanishing claim, note that there is no associated prime of $\tilde{\Gamma}'$ over $\iota(Z)$ and that $\iota(Z)$ is a Cartier divisor. For the second, recall that we showed earlier that $(\tilde{m}')^{-1}(\iota(Z))$ is smooth so the result follows from our $\TorTor$ vanishing result (proposition \ref{BertiniTor}).
\end{proof}

\section{Examples} \label{Examples}

In this section we will compute $g_M$ for several matroids $M$. First
of all, we observe that by our previous results we need only consider matroids which are not
direct sums, two-sums or series-parallel extensions of smaller
matroids. Every matroid can be built from these operations out of
three-connected matroids. (And, in a certain sense, uniquely so -- see \cite{CE}.) Therefore, in this section we will
only discuss computing $g_M$ for three-connected matroids.

\begin{figure}
$$\begin{matrix}
\Unif(2,4) & 2t+t^2 \\
\includegraphics[width=1 in]{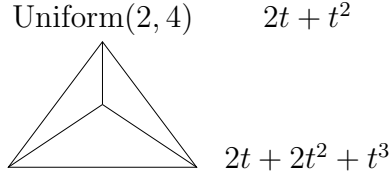} & 2t+2t^2+t^3 
\end{matrix}$$
\caption{The $3$-connected matroids with $\beta=2$} \label{beta2}
\end{figure}

\begin{figure}
$$\begin{matrix}
\Unif(2,5) & 3t+2t^2 \\
\includegraphics[width=1 in]{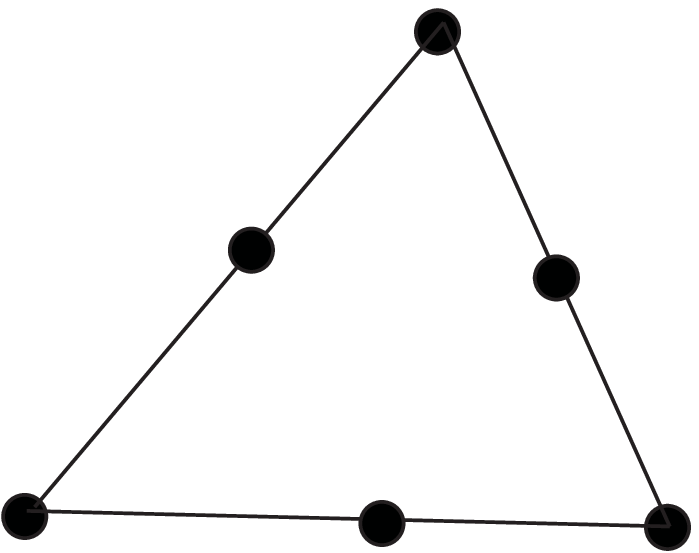} & 3t+3t^2+t^3 \\ 
\includegraphics[width=1 in]{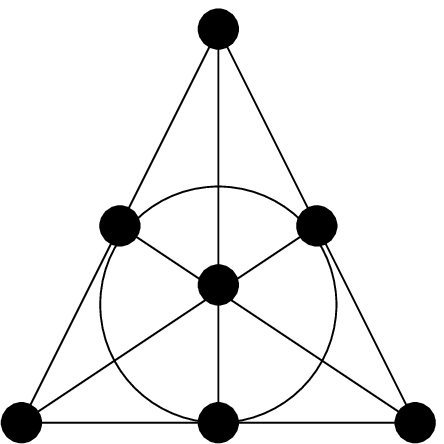} & 3t+5t^2+3t^3 \\
\includegraphics[width=1 in]{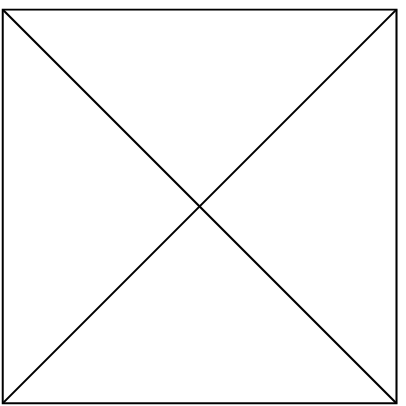} & 3t+5t^2+4t^3+t^4 
\end{matrix}$$
\caption{The $3$-connected matroids with $\beta=3$} \label{beta3}
\end{figure}

\begin{figure}
$$\begin{matrix}
\Unif(2,6) & 4t+3t^2 \\
\includegraphics[width=1 in]{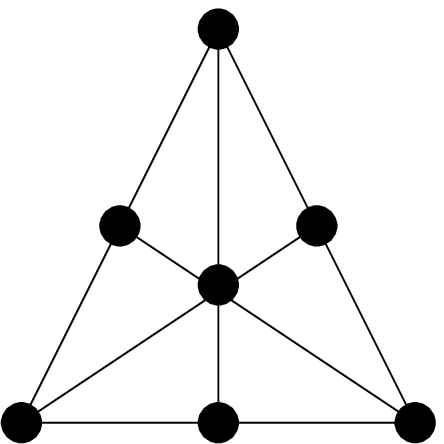} & 4t+6t^2+3t^3 \\
\includegraphics[width=1 in]{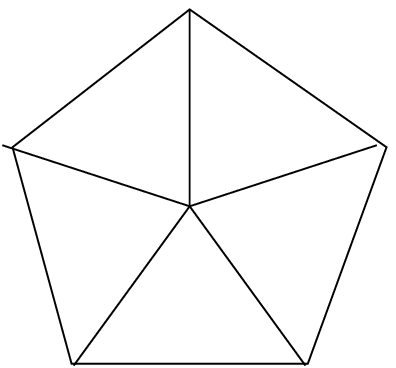} & 4t+9t^2+10 t^3+5 t^4+ t^5 \\
\begin{pmatrix} 1 & 0 & 0 & 0 & 1 & 0 & 0 & - \alpha \\
0 & 1 & 0 & 0 & -1 & 1 & 0 & 0 \\
0 & 0 & 1 & 0 & 0 & -1 & 1 & 0 \\
0 & 0 & 0 & 1 & 0 & 0 & -1 & 1 \end{pmatrix}, \quad \alpha \neq 0,1 & 4t+6t^2+4t^3+t^4 \\
\includegraphics[width=1 in]{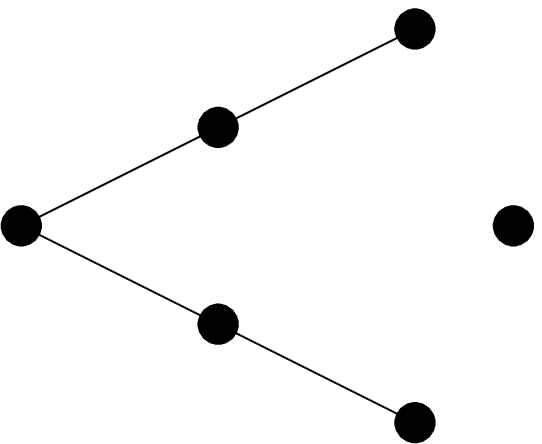} & 4t+4t^2+t^3 \\
\includegraphics[width=1 in]{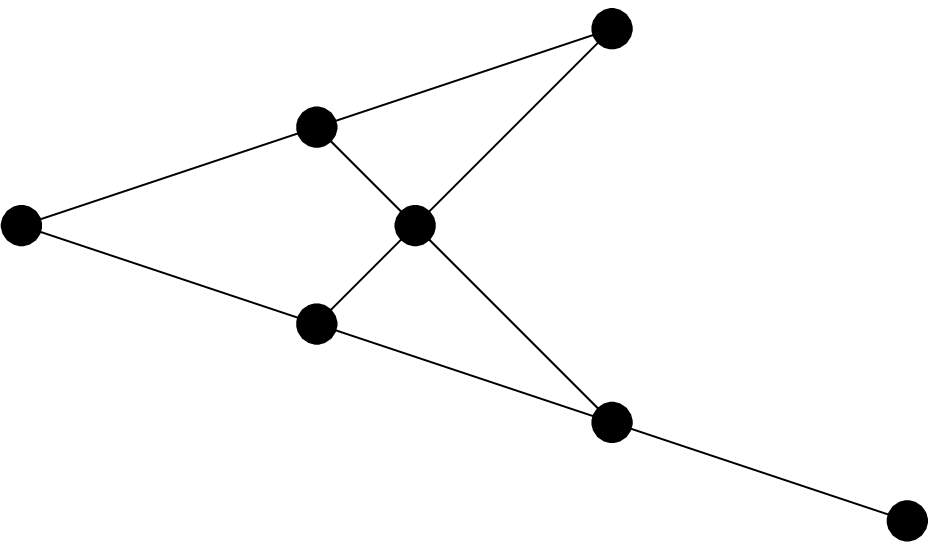} & 4t+5t^2+2t^3 \\
\includegraphics[width=1 in]{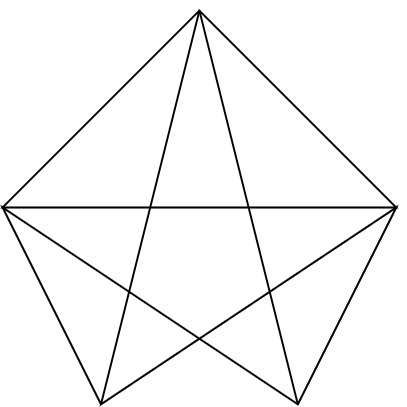} & 4t+8t^2+7t^3+2t^4 \\
\begin{pmatrix} 1 & 0 & 0 & 0 &1 & 1 & 0 & 1 \\
0 & 1 & 0 & 0 & 1 & 0 & 1 & 1 \\
0 & 0 & 1 & 0 & 0 & 1 & 1 & 1 \\
0 & 0 & 0 & 1 & 1 & 1 & 1 & 1 \\ \end{pmatrix}, \quad \textrm{char\ } k=2 & 4 t + 14 t^2 + 12 t^3+ t^4
\end{matrix}$$
\caption{The $3$-connected matroids with $\beta=4$} \label{beta4}
\end{figure}

Oxley has shown (see \cite{Ox}) that there are
only finitely many three-connected matroids with given
$\beta$-invariant and has enumerated those with $\beta$-invariant less
than or equal to four. In tables \ref{beta2}, \ref{beta3} and \ref{beta4}, we list $g_M$ for each
matroid in Oxley's list. In the first column, we list the matroid $M$. Our notation is as follows: if $M$ is a graphical matroid, we give a graph that represents it. (See \cite{White1}, chapter 6.) We denote by $\Unif(d,n)$ the uniform matroid of rank $d$ on $n$ elements -- the matroid for which every $d$-element subset of $[n]$ is a basis. If $M$ is rank $3$, we give an arrangement of points in the plane that represents $M$. (See \cite{White1}, section 1.1.A) To distinguish planar point arrangements from graphs, we place bold dots for the points in a planar point arrangement and not for the vertices of a graph. If $M$ can not be represented in any of these forms, we give a matrix whose row span has matroid $M$. (There are matroids which are neither graphical, uniform, rank $3$ nor realizable, but none of them appear in Oxley's list.) As $g_{M^{\perp}}=g_M$, we only list one of $M$ and $M^{\perp}$.

Most of the computations in these tables are consequences of results stated later in this section. Those that are not were carried out by finding a polyhedral subdivision of $\Delta(d,n)$ which contained the appropriate polytope as a facet and for which $g_M$ for all of the other faces could be computed more easily, often by recognizing them as two-sums. 

\textbf{Example:} Let $M$ be the Pappus matroid -- the rank 3 matroid on $[9]$ whose bases are all three element subsets of $[n]$ except for $123$, $456$, $789$, $159$, $168$, $249$, $267$, $348$ and $357$. Let $P_I$ be $0$ if $I$ is a basis of $M$ and $1$ otherwise. Then $P_I$ is a tropical pl\"ucker vector. The corresponding subdivision $\DD_P$ corresponds to taking $\Delta(3,9)$ and at the nine vertices corresponding to the non-bases, cutting off each of these vertices with all of its neighbors. One facet of $\DD_P$ corresponds to $M$, the other nine correspond to series-parallel matroids. There are $9$ internal faces of $\DD_P$ in codimension $1$ and these each correspond to direct sums of two series-parallel matroids. So $g_{M}=g_{\Unif(3,9)} - 9 t -9 t^2=(21 t + 30 t^2 + 10 t^3) - 9t - 9t^2=12 t+21 t^2+10 t^3$.


We now move to the (few) infinite classes of three-connected matroids
for which we can compute $g_M$. The first case we compute is one that
we need to prove our Main Theorem -- the case of a uniform matroid.

\begin{prop} \label{Uniform}
Let $\Unif(d,n)$ be the uniform matroid of rank $d$ on $n$ elements. (That is, the matroid of a generic element of $G(d,n)$.) Then
$$g_{\Unif(d,n)}=\sum_i \frac{(n-i-1)!}{(d-i)!(n-d-i)!(i-1)!} t^i.$$
\end{prop}

\begin{proof}
  The polytope associated to the uniform matroid is the hypersimplex
  $\Delta(d,n)$. In \cite{TLS}, I exhibit decompositions of
  $\Delta(d,n)$ into series-parallel matroidal polytopes with the
  required $f$-vector, thus proving this result. We prefer, however,
  to give a more geometric proof.
  
  Let $L \subset \CC^n$ be a $d$-plane corresponding to the uniform
  matroid; \emph{i.e.} all of the Pl\"ucker coordinates $p_I(L)$ are
  nonzero. Recall the rational map $m: \PP(L) \times \PP(L^{\perp})
  \dashrightarrow Z$. We claim that, in this case, the map is actually
  well defined. Suppose to the contrary that $(x_1: \cdots: x_n)
  \times (y_1: \cdots: y_n)$ is a point of $\PP(L) \times
  \PP(L^{\perp})$ where $m$ is not defined. This implies that, for all
  $i \in [n]$, $x_i y_i=0$. 
  
  Let $F \subset [n]$ be the set of $i$ for
  which $x_i=0$ and $G \subset [n]$ be the set of $i$ for which
  $y_i=0$, so $F \cup G=[n]$. Because all the Pl\"ucker coordinates of $L$ are nonzero, there is no point in $\PP(L)$ where $d$ of the coordinate functions simultaneously vanish.  Thus, $|F| \leq d-1$. Similarly, $|G| \leq n-d-1$. But then $|F| + |G| \leq n-2$, contradicting $[n]=F \cup G$. We conclude that $m : \PP(L) \times \PP(L^{\perp}) \to Z$ is well defined.
  
  Now, $\PP(L) \times \PP(L^{\perp})$ embeds in $\PP(L \times L^{\perp})$ by the Segre embedding, and the map $m$ is just the restriction to $\PP(L) \times \PP(L^{\perp})$ of a linear projection  $\Lambda: \PP(L \times L^{\perp}) \dashrightarrow Z$. Our goal is thus to compute the holomorphic Euler characteristic of the intersection with $\PP(L) \times \PP(L^{\perp})$ of the pull back from $Z$ of linear subspaces. Now $\Lambda$ is not a generic linear projection, because it has the property that the pull back to $\PP(L) \times \PP(L^{\perp})$ of the coordinate functions on $Z$ are reducible hypersurfaces, where as for a generic projection these would be smooth. However, if $\Lambda' : \PP(L \otimes L^{\perp}) \dashrightarrow Z$ is a generic  linear projection, we may find a flat family $\Lambda_{t}$  of linear projections joining $\Lambda$ to $\Lambda'$ and such that every projection in the family is well defined on $\PP(L) \times \PP(L^{\perp})$.  Let $m_t$ denote the restriction of $\Lambda_t$ to $\PP(L) \times \PP(L^{\perp})$, then $m_t^{-1}$ of a generic $\PP^{i-1}$ in $Z$ will be independent of $t$. So we are reduced to computing $\chi((m')^{-1}(\PP^{i-1}))$ for a generic $\PP^{i-1}$ in $Z$ and generic map $\Lambda': \PP(L \otimes L^{\perp}) \dashrightarrow Z$. In other words, we must compute the holomorphic Euler characteristic of the intersection of $\PP(L) \times \PP(L^{\perp})$ with a generic $d(n-d)-n+i$ plane. Let $H$ denote this generic $d(n-d)-n+i$ plane.

Now, the Segre embedding  $\PP(L) \times \PP(L^{\perp}) \into \PP(L \times L^{\perp})$ corresponds to the polytope $\Delta^{d-1} \times \Delta^{n-d-1}$ where $\Delta^i$ is the $i$-dimensional simplex. This polytope has a subdivision into unimodular simplices with $\frac{(n-j-1)!}{(d-j)!(n-d-j)! (j-1)!}$ interior faces of dimension $n-j-1$. This leads to a flat degeneration of $\PP(L) \times \PP(L^{\perp})$ to the union of coordinate subspaces which is the Stanley-Riesner complex associated to this subdivision. (See, for example, section 6 of \cite{GrobTor}.) Now, intersecting this flat family with $H$, we get a flat degeneration of $(\PP(L) \times \PP(L^{\perp}) ) \cap H$ to a linear subspace arrangement. This subspace arrangement has $\frac{(n-j-1)!}{(d-j)!(n-d-j)! (j-1)!}$ interior faces of dimension $i-j$ for $i \geq j \geq 1$ and the link of each interior faces is a sphere. Hence, the holomorphic Euler characteristic of this arrangement is $\sum_{j \leq i} (-1)^{j-1} \frac{(n-j-1)!}{(d-j)!(n-d-j)! (j-1)!}$ and, as holomorphic Euler characteristic is constant in proper flat families, 
$$\chi \left( H \cap (\PP(L) \times \PP(L^{\perp})) \right) = \sum_{j \leq i} (-1)^{j-1} \frac{(n-j-1)!}{(d-j)!(n-d-j)! (j-1)!}.$$

A little algebra now yields the desired claim.
\end{proof}

We next compute $g$ for the $d$-wheel and the $d$-whirl. The $d$-wheel, denoted $W_d$, is the graphical matroid associated to the graph consisting of a cycle of length $d$ and one additional vertex which is joined to every vertex in the cycle. It has $2d$ elements, rank $d$ and is represented by the matrix
$$\setcounter{MaxMatrixCols}{20}
M_d:=\begin{pmatrix}
1 & 0 & 0 & \cdots & 0 & 0 &      1 &  0 & 0 & \cdots & 0 & -1 \\
0 & 1 & 0 & \cdots & 0 & 0 &     -1 &  1 & 0 & \cdots & 0 & 0 \\
0 & 0 & 1 & \cdots & 0 & 0 &      0 & -1 & 1 & \cdots & 0 & 0 \\
\vdots & \vdots & \vdots & \ddots & \vdots & \vdots &  \vdots & \vdots & \vdots & \ddots & \vdots & \vdots \\
0 & 0 & 0 & \cdots & 1 & 0 & 0 & 0 & 0 & \cdots &  1 & 0 \\
0 & 0 & 0 & \cdots & 0 & 1 & 0 & 0 & 0 & \cdots & -1 & 1
\end{pmatrix}$$
The $d$-whirl, denoted $\mathcal{W}_d$, is the matroid on a $2d$
element set which has all of the bases that the $d$-wheel does and, in
addition, has the $d$ edges of the outer rim of the wheel as a basis. From this description it
is easy to see that $P_{\mathcal{W}_d}$ can be cut into two pieces,
one of which is $P_{W_d}$ and the other of which is $P_S$ where $S$ is a 
series-parallel matroid. These pieces meet along $P_T$, where $T$ is the direct sum of a two series-parallel matroids.  Thus, we see that $g_{\mathcal{W}_d}=g_{W_d}+t+t^2$. We will concentrate on computing $g_{\mathcal{W}_d}$.

The $d$-whirl is represented by the matrix
$$\setcounter{MaxMatrixCols}{20}
\mathcal{M}_d:=\begin{pmatrix}
1 & 0 & 0 & \cdots & 0 & 0 &      1 &  0 & 0 & \cdots & 0 & -\alpha \\
0 & 1 & 0 & \cdots & 0 & 0 &     -1 &  1 & 0 & \cdots & 0 & 0 \\
0 & 0 & 1 & \cdots & 0 & 0 &      0 & -1 & 1 & \cdots & 0 & 0 \\
\vdots & \vdots & \vdots & \ddots & \vdots & \vdots &  \vdots & \vdots & \vdots & \ddots & \vdots & \vdots \\
0 & 0 & 0 & \cdots & 1 & 0 & 0 & 0 & 0 & \cdots &  1 & 0 \\
0 & 0 & 0 & \cdots & 0 & 1 & 0 & 0 & 0 & \cdots & -1 & 1
\end{pmatrix}$$
where $\alpha$ is not $0$ or $1$.

\begin{prop}
$$g_{\mathcal{W}_d}(t) = (1+t)^d-1\quad \quad g_{W_d}(t) = (1+t)^d -1 -t-t^2$$
\end{prop}

\begin{proof}
As explained above, $g_{\mathcal{W}_d}=g_{W_d}+t+t^2$ so it is enough to prove the claimed equality for $g_{\mathcal{W}_d}$.  Let $L$ be the row span of $\mathcal{M}_d$. Then $L^{\perp}$ is the row span of the matrix
$$\setcounter{MaxMatrixCols}{20}
\mathcal{M}^{\perp}_d := \begin{pmatrix}
     1 &  -1 & 0 & \cdots & 0 & 0 & -1 & 0 & 0 & \cdots & 0 & 0   \\
  0 &  1 & -1 & \cdots & 0 & 0 & 0 & -1 & 0 & \cdots & 0 & 0    \\
     0 & 0 & 1 & \cdots & 0 & 0 & 0 & 0 & -1 & \cdots & 0 & 0   \\
\vdots & \vdots & \vdots & \ddots & \vdots & \vdots &  \vdots & \vdots & \vdots & \ddots & \vdots & \vdots \\
 0 & 0 & 0 & \cdots &  1 & -1 & 0 & 0 & 0 & \cdots & -1 & 0  \\
 -\alpha & 0 & 0 & \cdots & 0  & 1 & 0 & 0 & 0 & \cdots & 0 & -1
\end{pmatrix}$$

We coordinatize $L$ and $L^{\perp}$ via these matrices. The map $\PP(L) \times \PP(L^{\perp}) \to Z$ is then given by 
\begin{multline*}
m : (u_1 : \cdots : u_d) \times (v_1 : \cdots : v_d) \mapsto \\ (u_1v_1 -
\alpha u_1 v_d : u_2 v_2 - u_2 v_1 : \cdots : u_d v_d - u_d v_{d-1} : \\
u_2 v_1 - u_1 v_1 : u_3 v_2 - u_2 v_2 : \cdots : \alpha u_1 v_d - u_d
v_d ).\end{multline*}
We must compute the inverse image of an $(i-1)$-plane in $Z$. We
first consider the problem of computing the inverse image of the point
$(z_1: \cdots : z_{2d})$.

In this paragraph, we only work with maps up to birational isomorphisms. We can factor $m$ as the monomial map $\mu :\PP^{d-1} \times \PP^{d-1} \dashrightarrow \PP^{2d-1}$ given by 
$$\mu : (u_1 : \cdots : u_d) \times (v_1 : \cdots : v_d) \mapsto (u_1 v_1 :
\cdots : u_d v_d : \alpha u_1 v_d : u_2 v_1 : \cdots : u_d v_{d-1})$$
and the linear map $\Lambda : \PP^{2d-1} \dashrightarrow Z$ which is projection from
$(1:1:\cdots:1)$. Then $\Lambda^{-1}(z_1 : \cdots : z_{2d-1})$ is the closure of
the line of points of the form $(w_1+t : \cdots : w_{2d}+t) \in
\PP^{2d-1}$, where $w$ depends linearly on $z$. The point $(w_1+t : \cdots : w_{2d}+t) \in
\PP^{2d-1}$  is in the
closure of the image of $\mu$ if and only if $(w_1+t)(w_2+t) \cdots (w_d+t) = \alpha (w_{d+1}+t) (w_{d+2}+t) \cdots (w_{2d}+t)$. 

If $z$ runs over an $(i-1)$-plane, so does $w$. Thus the
$\tilde{m}$-preimage of this $(i-1)$-plane is birational to the
hypersurface $(w_1+t)(w_2+t) \cdots (w_d+t) = \alpha (w_{d+1}+t)
(w_{d+2}+t) \cdots (w_{2d}+t)$ in $\PP^{i}$, where each $w_r$ is a
generically chosen linear function of the $i$ coordinates other than
$t$.  One can check that, for $w_r$ chosen
generically, this is a smooth hypersurface of degree $d$ (here we use $\alpha \neq 1$). Therefore, we may use it to compute the
holomorphic Euler characteristic of $\tilde{m}^{-1}(\PP^{i-1})$.

A smooth hypersurface of degree $d$ in $\PP^{i}$ has holomorphic
Euler characteristic $1 - (-1)^i \binom{d-1}{i}$. Note that formula
gives $1$ for $i \geq d-1$, which is in accord with our conventions. So 
\begin{multline*}
\sum_{i=1}^{\infty} \chi(\tilde{m}^{-1}(\PP^{\min(n-2,i-1)})) t^i =\sum_{i=1}^{\infty} \left( 1-(-1)^i \binom{d-1}{i} \right) t^i= \\ \frac{t}{1-t} + \left( 1-(1-t)^{d-1}\right)  = \frac{1-(1-t)^d}{1-t} \end{multline*}
and $g_{\mathcal{W}_d}(t)=(1+t)^d-1$. 
\end{proof}

When $M$ has rank $2$ it is easy to describe $g_M$ : as $M$ is unaffected by parallel extension, we may assume that $M$ is not a parallel extension of any smaller matroid. We assume, as always, that $M$ contains no loops. Then $M$ is a uniform matroid, say on $n$ elements. Our computation above gives $g_M=(n-2)t+(n-3)t^2$.

We give a formula for $g_M$ in the case where $M$ has rank $3$, although we omit a detailed argument as overly lengthy. Since $g_M$ is unaltered by parallel extensions, we may assume that all the parallel classes of $M$ contain only one element. Also, we assume as always that $M$ contains no loops. Then $M$ has $n$ flats of rank $1$. Let $r$ be the number of flats of rank $2$ and let their cardinalities be $d_1$, \dots, $d_r$. (If $M$ is realizable, the rank two flats of $M$ are the vertices of the corresponding hyperplane arrangement in $\PP^2$ and the order of a flat is the number of hyperplanes passing through a vertex.) We have $\sum \binom{d_i}{2}=\binom{n}{2}$.

\begin{prop}
With the above assumptions and notations, we have 
\begin{multline*}
g_M=\left( \binom{n-2}{2} - \sum_i \binom{d_i-1}{2} \right) t +\left( (n-3)(n-4) - \sum_i (d_i-2)^2   \right) t^2 \\
+ \left( \binom{n-4}{2} - \sum \binom{d_i-2}{2} \right) t^3.
\end{multline*}
\end{prop}

\begin{proof}[Sketch of Proof]
The coefficient of $t$ is $\beta(M)$, which may be computed by any number of standard means.  We know that $g_M(-1)=-1$ by proposition \ref{ValueAtMinusOne}. We thus have two linear equations relating the three coefficients of $g_M$ and we will be able to determine $g_M$ as soon as we know one more linear relation between the coefficients. We set our sights on proving that the coefficient of $t^2$ in $(1-t) g_M$, in other words, $\chi(\tilde{m}^{-1}(\PP^1))$, is $1-\binom{n-4}{2}-\sum \binom{d_i-2}{2}$. 

Let $L$ be a $3$-plane in $n$-space corresponding to $M$, then we have a projection map $\tilde{\Gamma} \to \Gamma \to \PP(L) \times \PP(L^{\perp}) \to \PP(L)$. We will consider the image of $\tilde{m}^{-1}(\PP^1)$ in $\PP(L) \isomorph\PP^2$. This will be a curve $C$; we claim that $C$ has degree $n-3$ and its only singularities are ordinary multiple points. More specifically, there is a point of multiplicity $d_i-2$ for each $i$ such that $d_i \geq 4$ and no other singularities.

Roughly speaking, the argument is to apply standard elimination theory techniques to the projection $\PP(L) \times \PP(L^{\perp}) \to \PP(L)$ to get a determinantal formula for $C$  and show that $C$ obeys an equation of the form $\sum a_I \prod_{j \in [n] \setminus I} x_j=0$. Here $x_j$ is understood as the restriction to $\PP(L)$ of the coordinate function on $\PP^{n-1}$ and $I$ runs over the bases of $M$. We use the fact that $M$ has no parallel elements to show that this polynomial is irreducible and thus defines $C$. This polynomial vanishes to order $d_i-2$ at the point of $\PP(L)$ corresponding to the $i^{\textrm{th}}$ flat. A rather detailed computation is required to show that it has no further singularities. One must then check that $\tilde{m}^{-1}(\PP^1)$ is the normalization of $C$ and hence has holomorphic Euler characteristic $1-\binom{n-4}{2}-\sum \binom{d_i-2}{2}$. 
\end{proof}

\textbf{Remark:} It is also possible to prove this proposition by polyhedral combinatorics. Specifically, using our assumption that $M$ has no parallel elements, the function $p_I=0$ if $I$ is a basis of $M$ and $p_I=1$ if $I$ is not is a tropical Pl\"ucker vector. The corresponding matroidal subdivision of $\Delta(3,n)$ has $r+1$ facets: one facet corresponds to the matroid $M$ and the others correspond to parallel extensions of $\Unif(3,d_i+1)$. This subdivision has $r$ interior faces of codimension $1$, each corresponding to $\Unif(2,d_i) \oplus \Unif(1,n-d_i)$. There are no faces of codimension greater than $1$.

\section{Further Questions and Speculations}

The most obvious challenge is to extend the results of this paper to matroids realizable only in characteristic $p$ or not realizable at all. We conjecture that they can be so extended. Specifically, 

\begin{conj}
For any matroids $M$, $M_1$ and $M_2$, we have $g_{M}=g_{M^{\perp}}$, $g_{M_1 \oplus M_2}=g_{M_1} g_{M_2}$, $g_{M_1 +_2 M_2}=t g_{M_1} g_{M_2}$ and all the coefficients of $g_M$ are nonnegative.
\end{conj}

Let $\Flag(1,d,n-1;n)$ be the variety of partial flags of dimensions $(1,d,n-1)$ in $\CC^n$. Let $p$ and $q$ be the projections of $\Flag(1,d,n-1;n)$ to $G(d,n)=\Flag(d;n)$ and to $\Flag(1,n-1;n)$ respectively. Let $x \in G(d,n)$. Then the coefficient of $t$ in $g_x(t)$ is the degree of the map $p^{-1}(\overline{Tx}) \to \Flag(1,n-1;n)$ and the other coefficients of $g_x$ are computed in terms of holomorphic Euler characteristics of $p^{-1}(\overline{Tx}) \cap q^{-1}(\Omega)$ for various subvarieties $\Omega \subset \Flag(1,n-1;n)$. This shows that $g_x$ is determined by the class $q_* p^* [\OO_{\overline{T x}}]$ in $K^{\circ}(\Flag(1,n-1;n))$. Now, $\Flag(1,n-1;n)$ embeds as a hypersurface in $\PP^{n-1} \times \PP^{n-1}$; let $\iota: \Flag(1,n-1;n) \into \PP^{n-1} \times \PP^{n-1}$ denote this embedding. Then $\iota_* : K^{\circ}(\Flag(1,n-1;n)) \to K^{\circ}(\PP^{n-1} \times \PP^{n-1}) = \ZZ[s,t]/{<s^n, t^n>}$ is injective. Here $s$ and $t$ are the pullbacks of the hyperplane classes from the first and second factor of $\PP^{n-1} \times \PP^{n-1}$. 

Let $\gamma_x(s,t) \in \ZZ[s,t]/{<s^n, t^n>}$ represent $\iota_* q_* p^*  [\OO_{\overline{T x}}]$. It is easy to check that $(-1)^c \gamma_x(0,-t)=g_x(t)$, so $\gamma$ determines $g$. From this perspective, $\gamma$ seems like a more natural object than $g$. However, it appears that $\gamma$ actually contains no additional information. 
\begin{conj}
We have $\gamma_x(s,t)= (-1)^c g_x(-s-t+st)$. In particular, $\gamma_x(s,t)=\gamma_x(t,s)$.
\end{conj}
In particular, since it is easy to show that $\gamma_{x}(s,t)=\gamma_{x^{\perp}}(t,s)$, this would yield another proof that $g_{x}(t)=g_{x^{\perp}}(t)$. 

We would like to comment on the plausible geometric basis of this conjecture. We saw $g_x(t)$ can be computed by finding $\chi(\tilde{m}^{-1}(\PP^{i-1}))$ for a generic $\PP^{i-1}$. Similarly, $\gamma_x(s,t)$ could be computed if we knew $\chi(\tilde{m}^{-1}(\PP^{i-1} \times \PP^{j-1}))$; here the embedding $\PP^{i-1} \times \PP^{j-1} \into Z$ is given by $(u_1 : \ldots : u_i) \times (v_1: \ldots v_j) \mapsto (a_1(u) b_1(u) : \ldots : a_n(u) b_n(v))$ where $a_k$ and $b_k$ are linear functions chosen generically with respect to the constraint $\sum a_k(u) b_k(v)=0$. We could instead take an embedding $\PP^{i-1} \times \PP^{j-1} \into Z$ by $(u_1 : \ldots : u_i) \times (v_1: \ldots v_j) \mapsto (q_1(u,v) : \ldots : q_n(u,v))$ where the $q_k$ are bilinear forms chosen generically subject to the requirement $\sum q_k=0$. With this modified embedding it is possible to compute $\chi(\tilde{m}^{-1}(\PP^{i-1} \times \PP^{j-1}))$  and get a result that would imply the conjecture. The difficulty is to show that the requirement that $q_k(u,v)$ factor as $a_k(u) b_k(v)$ does not alter the computation.

In this paper we have presented a number of properties of $g_M$ and it is natural to ask whether they are enough to uniquely determine $g_M$. Experimentation suggests that even a small subset of these properties suffice. Specifically, we make the following conjecture, which is essentially combinatorial:

\begin{conj}
Let $M \mapsto G_M$ be a map from isomorphism classes of loop and co-loop free matroids to $\ZZ[t]$ obeying the following properties:
\begin{enumerate}
\item If $\DD$ is a polyhedral decomposition of the matroid polytope $\Poly_M$ into smaller matroid polytopes then $G_M=\sum_{\Poly_F \in \oDD} G_F$ where $\oDD$ is the set of interior faces of $\DD$.
\item If $M$ is a direct sum of $c$ series-parallel matroids then $G_M=t^c$.
\end{enumerate}
Then $G_M=g_M$.
\end{conj}

One major flaw with our results so far is that we have not presented a simple recursion for $g_M$ that would allow it to be efficiently computed. In searching for such a recursion, I attempted to modify the argument used to prove proposition \ref{Additive}. Recall the notations $x/e$ and $x \setminus e$ introduced in section \ref{BetaInvariant}; at times we will attach the subscripts $x$, $x/e$ or $x \setminus e$ to an object or map to indicate which of these points it is associated to.  My strategy was as follows: let $H$ be a projective $(i-1)$-space in $Z$ which is generic subject to the condition that $z_e$ is zero on $H$.  One must then show that $\chi(\tilde{m}^{-1}(H))$ is the same as $\chi(\tilde{m}^{-1}(\PP^{i-1}))$ for a generic $\PP^{i-1}$ in $Z$. Assuming this holds, the next question is to figure out what $\tilde{m}^{-1}(\{ z_e=0 \})$ is after which we can hope to determine $\chi(\tilde{m}^{-1}(H))$ for a generic $H \subset \{ z_e=0 \}$. 

What appears to happen is that $\tilde{m}^{-1}(\{ z_e =0 \})$ has two components, one isomorphic to $\tilde{\Gamma}_{x/e}$ and the other to $\tilde{\Gamma}_{x \setminus e}$. Let us assume that this is correct for the remainder of this section. When $i=1$, we can deduce that $\tilde{m}_x^{-1}(H) = \tilde{m}_{x/e}^{-1}(H) \sqcup \tilde{m}_{x \setminus e}^{-1}(H)$ and thus deduce proposition \ref{Additive}. Once $i$ is larger than $1$, we will still have $\tilde{m}_x^{-1}(H) = \tilde{m}_{x/e}^{-1}(H) \cup \tilde{m}_{x \setminus e}^{-1}(H)$ but the union will not be disjoint. Thus, we must understand how $\tilde{\Gamma}_{x/e}$ and $\tilde{\Gamma}_{x \setminus e}$ meet inside $\tilde{m}_x^{-1}( \{ z_e =0 \})$.

Suppose that $L_1 := L(x_1) \subseteq L_2 := L(x_2)$ are two subspaces of $\CC^n$, with $L_1$ not containined in a coordinate subspace and $L_2$ not containing a coordinate axis. Then we have the rational map $m : \PP(L_1) \times \PP(L_2^{\perp}) \dashrightarrow Z$ defined as before and we can resolve the singularities of the graph of $m$ to produce a a variety $\tilde{\Gamma}_{x_1,x_2}$ with a map to $Z$. In this vocabulary, what appears to occur is that $\tilde{m}_x^{-1}(\{ z_e =0 \})$ consists of two components, isomorphic to  $\tilde{\Gamma}_{x/e}$ and $\tilde{\Gamma}_{x \setminus e}$ and glued along $\tilde{\Gamma}_{x/e, x \setminus e}$. If this is true, then we would have $g_x=g_{x/e}+g_{x \setminus e} + g_{x/e, x \setminus e}$. where $g_{x_1,x_2}$ is the obvious generalization of $g_x$. 

This raises the hope of finding a recursive formula for $g$ which involves not only the Grassmannian but also two step flag manifolds. (A note for combinatorialists -- the matroid analogue of a two step flag is a pair of matroids related by a strong map.) The most natural guess would be that the preimage of $\{ z_e=0 \}$ in  $\tilde{\Gamma}_{L_1,L_2}$ would consist of a copy of  $\tilde{\Gamma}_{L_1/e,L_2/e}$ and a copy of  $\tilde{\Gamma}_{L_1 \setminus e,L_2 \setminus e}$ glued along $\tilde{\Gamma}_{L_1/e,L_2 \setminus e}$. Each of these components do in fact occur, but, in general, so do other components.  Figuring out a combinatorial description of the other components that occur seems to be the major obstacle to pursuing this idea towards an efficient recurrence. 

\textbf{Example:} Let us try to use this strategy to compute $g_{W_4}$  for $W_4$ the $4$-wheel. Chooising $e$ to be one of the edges in the rim of the wheel, $W_4 \setminus e$ is series-parallel and $W_4 / e$ is a parallel extension of $W_3$. So we have $g_{W_4}=(2t+2t^2+t^3)+t+g_{L_1, L_2}$ where here $L_1$ and $L_2$ are the row-spans of the following matrices:
$$L_1 = \Span \begin{pmatrix} 1 & -1 & 0 & -1 & 0 & 1 & -1 \\ 0 & 1 & 1 & -1 & 0 & 0 & 0 \\ 0 & 0 & 0 & 1 & 1  & -1&  0 \end{pmatrix} \quad  L_2 = \Span \begin{pmatrix} 1 & -1 & 0 & -1 & 0 & 1 & -1 \\ 0 & 1 & 1 & -1 & 0 & 0 & 0 \\ 0 & 0 & 0 & 1 & 1  & -1 & 0 \\ 1 & -1 & 0 & 0 & 0 & 0 & 0 \end{pmatrix}.$$
Then we have 
$$L_2^{\perp} =\Span \begin{pmatrix}  0 & 0 & 1 & 1 & -1 & 0 & 0 \\ 0 & 0  & 0 & 0 & 1 & 1 & 1\\ 1 & 1 & -1 & 0 & 0 & 0 & 0  \end{pmatrix}.$$
The rational map $m$ is generically one to one and its image is the hypersurface 
$$(z_1)(z_6+z_7) (z_1+z_2+z_3) +(z_7)(z_1+z_2)(z_5+z_6+z_7)=0.$$
(Recall that we also have $z_1+z_2+\cdots+z_7=0$.)

This is the toric four-fold associated to the polytope  
gotten by taking a product of two triangles and deleting three pairwise nonadjacent vertices. Also, $\PP(L_1) \times \PP(L_2^{\perp}) \isomorph \PP^2 \times \PP^2$ is the toric variety associated to the product of two triangles and the rational map $m$ is the one corresponding to deleting the vertices. Let $\tilde{\Gamma}$ be the blow up of $\PP(L_1) \times \PP(L_2^{\perp})$ at the three points corresponding to the deleted vertices -- explicitly, these points are	 $(1,0,0) \times (1,0,0)$, $(0,1,0) \times (0,1,0)$, and $(0,0,1) \times (0,0,1)$. One can check that $m$ now extends to a well defined morphism $\tilde{m} : \tilde{\Gamma} \to Z$. 

Then $\tilde{m}^{-1}(\{ z_1=0 \})$ has three components: the proper transform of $u_1=0$, the proper transform of $v_3=0$ and the exceptional fiber resulting from the blow up of $(0,1,0) \times (0,1,0)$. Each of these is a rational three-fold. The first two of these components are $\tilde{\Gamma}_{x_1/e_1, x_2/e_1}$ and $\tilde{\Gamma}_{x_1 \setminus e_1, x_2 \setminus e_1}$ and their intersection is  $\tilde{\Gamma}_{x_1 / e_1, x_2 \setminus e_1}$. However, I know of no general combinatorial rule which would have predicted the third component.

These three-folds meet each other transversely and all three of them meet together in a $\PP^1$. The images of these three-folds under $\tilde{m}$ are the $\PP^{3}$'s cut out of $z_1=z_1+z_2+\cdots+z_7=0$ by the further relations $z_7=0$, $z_1+z_2=0$ and $z_5+z_6+z_7=0$ respectively. On each component of $\tilde{m}^{-1}(\{ z_1=0 \})$, and on each intersection of components, the map $\tilde{m}$ is a birational isomorphism onto its image. This shows us that a generic line in $\{ z_1=0 \}$ has three preimages, that the preimage of a general plane in $\{ z_1 = 0 \}$ is three $\PP^1$'s arranged in a ring and the preimage of $H$ for $\dim H \geq 3$ has holomorphic Euler characteristic $1$. Hence $g_{x_1,x_2}=(1+t)(3t^2-0t^3+t^4-t^5+t^6 - \cdots)=3t^2+3t^3+t^4$. This gives $g_{W_4}=(2t+2t^2+t^3) + t + (3t^2+3t^3+t^4)=3t+5t^2+4t^3+t^4$, which is correct. If we had not realized that the exceptional fiber was there, we would have thought that $g_{x_1,x_2}=2t^2+2t^3+t^4$ and gotten the wrong answer. This ends our example.

A question which I have often been asked when speaking on this material is whether $g_x$ is determined by the Tutte polynomial $t_{M(x)}(z,w)$ of $M(x)$. (See chapter 6 of \cite{White2} for background on the Tutte polynomial.) This is a particularly interesting question because lemma 6.4 of \cite{TLS} shows that the Tutte polynomial behaves in the same manner as $(-1)^c g_x$ in matroidal subdivisions of $\Poly_M$ and the Tutte polynomial is multiplicative in direct sums. Nevertheless, $(-1)^c g_x$ is not any linear function of the Tutte polynomial. This is true because, in section 6 of \cite{TLS}, I present three matroids $M_1$, $M_2$ and $M_3$ such that $t_{\Unif(3,6)}=6 t_{M_1} - 9 t_{M_2} + 4 t_{M_3}$ and the corresponding relation does not hold between the polynomials $(-1)^c g_x$. I very much suspect that $g_x$ is not determined by the Tutte polynomial, but I do not have an explicit example of this because it is difficult to find two nonisomorphic matroids with the same Tutte polynomial. This does, however, raise the following question, to which I would be curious to know the answer:

\textbf{Question:} Is there a linear map $K^{\circ}(G(d,n)) \to \ZZ[z,w]$ taking $[\OO_{\overline{Tx}}]$ to $t_{M(x)}(z,w)$? Does this map have any geometric significance?

\section{Appendix: Proofs of some toric lemmas} \label{Appendix}

In this section, we will prove some of the claims about toric varieties made earlier in the text. None of the material in this section is original, but some of it is hard to find in the published literature.

\begin{prop} \label{ProjNormal}
Let $x \in G(d,n)$. Then $\overline{Tx}$ is projectively normal. 
\end{prop}

\begin{proof}
Let $A$ be the subalgebra of $\CC[t_1,\ldots, t_n]$ generated by $p_{i_1 \cdots i_d}(x) t_{i_1} \cdots t_{i_d}$ for $(i_1, \ldots, i_d)$ ranging over $\binom{[n]}{d}$. Then $\overline{Tx}$ is $\Proj A$; we must show that $A$ is integrally closed. Clearly, $A$ is also the subalgebra generated by $t_{i_1} \cdots t_{i_d}$ where now $(i_1,\cdots, i_d)$ ranges over those elements of $\binom{[n]}{d}$ for which $p_{i_1 \cdots i_d}(x)$ is not zero -- in other words the bases of $M(x)$. This is the semigroup ring of the semigroup $S$ generated by $e_{i_1} + \cdots + e_{i_d}$, where $(i_1, \ldots, i_d)$ ranges over the bases of $M$. A semigroup ring is integrally closed if and only if the semigroup is saturated; that $S$ is saturated is shown in \cite{White3}.
\end{proof}

\begin{prop} \label{ToricFiber}
Let $R$, $K$ and $v$ be as in section \ref{Intro}. Let $x \in G(d,n)(K)$ with all the Pl\"ucker coordinates $p_I(x) \neq 0$. Let $P_I=v(p_I)$ and let $\DD_P$ be as in section \ref{Intro}. Let $\mathcal{X} \in G(d,n) \times \Spec R$ be the closure of $x$ and let $Y$ be the fiber of $\mathcal{X}$ over $\Spec \CC$. Then $Y$ is a union of toric varieties, indexed by and glued according to the faces of $\DD_P$.
\end{prop}

\begin{proof}
Let $Z$ denote the union of toric varieties glued along the faces of $\DD_P$.  Section 2 of \cite{Smirnov} shows that $Z$ is the radical of $Y$ (this is true for any regular subdivision of any lattice polytope). To check that the equality is one of schemes, we check that both objects have the same Hilbert function; this is enough because the complex of toric varieties is reduced. Now, $Y$ is a flat degeneration of $\overline{T x_{\gen}}$ where $x_{\gen}$ is a generic point of $G(d,n)$, so the Hilbert function of $Y$ is the same as that of $\overline{T x_{\gen}}$. Specifically, $h_Y(N)$ is the number of lattice points of the form $a_1+\cdots+a_N$ where each $a_i$ is a vertex of $\Delta(d,n)$ -- call the set of lattice points of this form $A_N$. On the other hand, $h_Z(N)$ is the number of lattice points of the form $a_1+\cdots +a_N$ where there is some particular face $F$ of $\Delta(d,n)$ such that each $a_i$ is a vertex of $F$ -- call the set of lattice points of this form $B_N$.

Clearly, $B_N \subseteq A_N$. Suppose now that $a \in A_N$. Then $a/N \in \Delta(d,n)$ and, in particular, lies in some face of $\DD_P$, say $F$. Then $a$ is a lattice point, the sum of whose coordinates is $Nd$, in the real cone spanned by $F$. But, by the result of \cite{White3}, the semigroup of lattice points in this cone whose coordinate sum is divisible by $d$ is generated by the vertices of $F$. So $A_N \subseteq B_N$, $A_N=B_N$ and $h_Y(N)=h_Z(N)$.
\end{proof}

\textbf{Remark:} Those readers who prefer algebraic arguments to geometric ones may like to read the extremely clear paper \cite{GrobTor}, which establishes the result of \cite{Smirnov} in the case that the regular subdivision involved is a triangulation. The proof of theorem 6.1 in that paper was the model for our argument here showing that the equality is one of schemes and not simply of point sets in our setting.

Let $Y$ be as above. Let $f_c$ be the number of $(n-c)$-dimensional interior faces of $\DD_P$ and, for each such face, let $x_j^c$ be an element of $G(d,n)$ such that $\overline{T x_j^c}$ is the stratum of $Y$ corresponding to the face. 

\begin{prop} \label{ToricComplex}
The complex of sheaves
$$0 \to \OO_Y \to \bigoplus_{j=1}^{f_1} \OO_{\overline{T x_j^1}} \to \cdots \to \bigoplus_{j=1}^{f_c} \OO_{\overline{T x_j^c}} \to \cdots \to \bigoplus_{j=1}^{f_n} \OO_{\overline{T x_j^n}} \to 0$$
is exact.
\end{prop}

\begin{proof}
Each of the schemes involved is $\Proj$ of a semigroup ring or, in the case of $Y$, $\Proj$ of a semigroup ring modulo a radical monomial ideal. We claim that this exactness holds even on the level of graded rings. This exactness is condition (3) of Theorem 4.2 in \cite{Ezra}; condition (1) of that theorem applies because $\DD_P$ is a subdivision of a ball.
\end{proof}

\begin{prop} \label{RationalSing}
For any $x \in G(d,n)$, $\overline{T x}$ has rational singularities.
\end{prop}

\begin{proof}
By proposition \ref{ProjNormal}, $\overline{Tx}$ is normal. By definition, $T$ acts on $\overline{T x}$ with a dense orbit. So we have checked that $\overline{T x}$ is a toric variety in the sense of \cite{Fult}. By section 2.6 of \cite{Fult}, this implies that $\overline{Tx}$ has rational singularities.
\end{proof}

\begin{prop} \label{MatroidOnly}
Let $x \in G(d,n)$ and let $M(x)$ be the matroid on $[n]$ whose bases are those $I$ for which $p_I(x) \neq 0$. Then the class of $\OO_{\overline{Tx}}$ in $K^{\circ}(G(d,n))$ can be determined from the isomorphism class of $M(x)$.
\end{prop}

To prove this, we will use the $T$-equivariant $K$-theory of $G(d,n)$. The standard reference for this subject is \cite{KK}. However, \cite{KK} does not provide the computational perspective we will need in this section. It seems to be difficult to find the results we need spelled out in a single place, so we summarize them in the following two paragraphs. None of the material in these paragraphs is original; in order to ease reading, references are consigned to the footnotes.

Let $\Lambda$ be the ring of degree zero Laurent polynomials in $n$ variables. The equivariant $K$-theory of $G(d,n)$, denoted $K^T(G(d,n))$, can be described as the ring of functions $I \mapsto f_I$ from $\binom{[n]}{d}$ to $\Lambda$ obeying the following condition: for every $B \in \binom{[n]}{d-1}$ and $i$, $j \in [n] \setminus B$, we have $f_{B \cup \{i \}} \equiv f_{B \cup \{ j \}} \mod 1-x_i/x_j$.  \footnote{This statement in the topological category after tensoring with $\QQ$ is a consequence of corollary A.5 of \cite{KR}; to get this result in the algebraic category without tensoring with $\QQ$ see corollary 5.12 of \cite{VV}. The reader should be warned that \cite{VV} works with the  ``higher" $K$ theory, which is far more subtle than $K^{\circ}$ but contains $K^{\circ}$ as its degree zero part.} Let $U_I$ be the open subset of $G(d,n)$ where $p_I \neq 0$.   It is standard that $U_I$ is a $d(n-d)$-dimensional affine space whose coordinates are naturally indexed by pairs $(i,j) \in I \times [n] \setminus I$. Furthermore, $U_I$ is {$T$-invariant} and $T$ acts on the $(i,j)$ coordinate by $t_j/t_i$. If $Z$ is any $T$-invariant subscheme of $G(d,n)$ then $Z \cap U_I$ is also $T$-invariant and thus has a $\ZZ^n$ graded Hilbert series $h_I(x_1, \ldots, x_n)$. The map $I \mapsto f_I$ corresponding to $[ \OO_Z ]$ is given by $f_I=h_I \prod_{i \in I} \prod_{j \in [n] \setminus I} \left( 1-x_j/x_i \right)$. \footnote{This result, including the fact that this formula gives a Laurent polynomial, is a special case of the general results in Section 8.2 of \cite{MSturm} -- see also the remarks at the end of chapter 8 of \cite{MSturm}.}

The ordinary $K$-theory of $G(d,n)$ is obtained from $K^T(G(d,n))$ by a canonical isomorphism $K^{\circ}(G(d,n)) \isomorph K^T(G(d,n)) \otimes_{\Lambda} \ZZ$, where every monomial of $\Lambda$ acts on $\ZZ$ by $1$. \footnote{This result in the topological category after tensoring with $\QQ$ is an easy consequence of the main result of \cite{KR}, combined with the analogous result in cohomology. Once again, to get into the algebraic category and remove the tensor product with $\QQ$, the reader should see Theorem 5.19 of \cite{VV}.} In particular, the class of $\OO_{\overline{Tx}}$ in $K^{\circ}(G(d,n))$ is determined by its class in $K^T(G(d,n))$.

We introduce the following notation: For $(a_1, \ldots, a_n) \in \ZZ^n$, we write $x^a$ for $x_1^{a_1} \ldots  x_n^{a_n}$. For $P$ a lattice polytope in $\RR^n$, we write $h_P$ for $\sum_{a \in P \cap \ZZ^n} x^a$. For $P$ a polytope in $\RR^n$ and $v \in \RR^n$, we write $v+P$ for the translation of $P$ by $v$.

\begin{proof}[Proof of proposition \ref{MatroidOnly}]
We first show that the class of $\OO_{\overline{Tx}}$ in $K^T(G(d,n))$ is determined by the data of $M(x)$ \textbf{and the labeling of the elements of $M(x)$ by $[n]$.}

Let $I$ be a basis of $M(x)$, let $S_I \subset \ZZ^n$ be the semigroup generated by the set of vectors of the form $e_j-e_i$ where $(i,j) \in I \times [n] \setminus I$ and $I \cup \{ j \} \setminus \{ i \}$ is also a basis of $M(x)$. Then $\overline{T x} \cap U_I$ is isomorphic to $\Spec \CC[S_I]$ and the Hilbert series of $\overline{T x} \cap U_I$ is $\sum_{a \in S_I} x^a$. We see that $f_I$ is completely determined by $M(x)$ and its labeling by $[n]$. If $I$ is not a basis of $M(x)$, then $f_I=0$. We see that the class of $\OO_{\overline{Tx}}$ in $K^{T}(G(d,n))$ is determined by $M(x)$ and its labeling by $[n]$.

Now, the class of $\OO_{\overline{Tx}}$ in $K^{\circ}(G(d,n))$ is determined by its class in $K^T(G(d,n))$.   We must show that the class of $\OO_{\overline{Tx}}$ in $K^{\circ}(G(d,n))$ is determined purely by the isomorphism class of $M(x)$ and not by the labeling of its elements by $[n]$. Suppose we labeled the elements of $M(x)$ differently. This would have the same effect as acting on $\overline{Tx}$ by a permutation matrix. But $GL_n$, and hence its subgroup $S_n$, acts trivially on $K^{\circ}(G(d,n))$, so this would have no effect on the class of $\OO_{\overline{Tx}}$ in $K^{\circ}(G(d,n))$.

\end{proof}

\textbf{Remark:} It follows from the proof of proposition \ref{MatroidOnly} that equation (\ref{KTheoryEquality}) holds for nonregular matroidal decompositions of $\Delta(d,n)$.

The following result shows that it makes sense to talk about $g_M$ for a nonrealizable matroid. We will not need this fact in this paper, but it makes the study of $g_M$ potentially far more interesting.

\begin{prop} \label{NonRealizable}
Let $M$ be a rank $d$ matroid on $n$ elements, which may or may not be realizable. Define a function $I \mapsto f_I$ from $\binom{[n]}{d}$ to $\Lambda$ by the recipe in the previous proof. Then $I \mapsto f_I$ obeys  $f_{B \cup \{i \}} \equiv f_{B \cup \{ j \}} \mod 1-x_i/x_j$ and thus defines a class in $K^{T}(G(d,n))$ and in $K^{\circ}(G(d,n))$.
\end{prop}

We will need several polyhedral lemmas before proving this result.

\begin{prop} \label{RaysAreDenominator}
Let $C$ be a (polyhedral, generated by lattice vectors) cone in $\ZZ^n$ which does not contain any line. Then $h_C$  is the quotient of a Laurent polynomial by $\prod_{b \in R(C)} (1-x^b)$ where $R(C)$ is the set of minimal lattice vectors along extremal rays of $C$.
\end{prop}

\begin{proof}
This is theorem 4.6.11 of \cite{Stan}, combined with the equality $\textbf{i} = \textbf{iii}$ in the preceding proposition 4.6.10. (Stanley includes the additional assumption that $C \subseteq \RR_{\geq 0}^n$, but it is easy to check that removing this assumption just changes the numerator from an ordinary polynomial to a Laurent polynomial.)
\end{proof}

\begin{prop} \label{PolytopeRaysAreDenominator}
Let $P$ be a lattice polytope. Then $h_P$ is a quotient of a Laurent polynomial by $\prod_{b \in R(P)} (1-x^b)$ where $R(P)$ is the set of minimal lattice vectors along the unbounded rays of $P$.
\end{prop}

\begin{proof}
Let $C \subset \ZZ^{n+1}$ be the cone on $P \times \{ 1 \}$. Then $R(C) = ( R(P) \times \{ 0 \}) \sqcup ( V(P) \times \{ 1 \} )$ where $V(P)$ is the set of vertices of $P$. By proposition \ref{RaysAreDenominator}, 
$$h_C(x_1, \ldots, x_n, x_{n+1}) = \frac{L(x_1^{\pm 1}, \ldots, x_n^{\pm 1}, x_{n+1}^{\pm 1})}{\prod_{b \in V(P)} (1- x^b x_{n+1}) \prod_{b \in R(P)} (1- x^b) }$$
where $L$ is a Laurent polynomial. (In fact, $L$ will only have nonnegative powers of $x_{n+1}$, but we will not need that.) Expanding $h_C$ as a power series in $x_{n+1}$, we see that the coefficient of any given power of $x_{n+1}$ has the required form. $h_P$ is the coefficient of $x_{n+1}$.
\end{proof}

For $P$ any lattice polytope and $v$ a vertex of $P$, let $T_v(P) \subset \RR^n$ denote the cone $\{ w : v+\epsilon w \in P \hbox{\ for\ $\epsilon>0$\ sufficiently\ small} \}$. The following is a result of Brion:

\begin{prop} \label{BrionsLemma}
Let $P$ be a lattice polytope in $\ZZ^n$ which does not contain any line. Then $h_P=\sum_{v \in V(P)} h_{v+T_v(P)}$.
\end{prop}

\begin{proof}
In the case where $P$ is bounded, this is proved very elegantly by methods of equivariant $K$-theory in \cite{Brion2}. We will need the result for unbounded $P$, however, for which we cite \cite{Ishida}.
\end{proof}

\begin{lemma} \label{EdgeCancellation}
Let $P$ be any polytope and let $u$ and $v$ be vertices of $P$ connected by an edge of $P$. Let $e$ be the minimal lattice vector along the edge pointing from $u$ to $v$, with $v=u+ke$. Then $h_{T_u(P)}+h_{T_v(P)}$ is a rational function whose denominator is \textbf{not} divisible by $1-x^e$.
\end{lemma}

\begin{proof}
We first prove the lemma under the following hypothesis: $(u+T_u(P)) \cap (v+T_v(P))$ is a polytope with its only vertices at $u$ and $v$. We will say $P$ is \emph{clean} if this property is satisfied. In this case, $h_{u+T_u(P)}+h_{v+T_v(P)}$ is, by Brion's result (Lemma \ref{BrionsLemma}), a rational function whose denominator is \emph{not} divisible by $1-x^e$.
Then
$$h_{T_u(P)}+h_{T_v(P)} = x^{-v} \left( (x^{ek}-1) h_{u+T_u(P)} + (h_{u+T_u(P)}+h_{v+T_v(P)}) \right).$$
Since $x^e-1$ only divides the denominator of $h_{u+T_u(P)} $ once, it does not divide the denominator of $ (x^{ek}-1) h_{u+T_u(P)}$. So we have shown the claim when $P$ is clean

We now consider the general case when $P$ is not clean. Let $y$ be a point in the interior of the edge $uv$ and let $C$ be the cone 
$\{ w : y+\epsilon w \in P \hbox{\ for\ $\epsilon>0$\ sufficiently\ small} \}$; note that $C$ contains $T_v(P)$ and $T_u(P)$. Now, choose a halfspace $H_v$ such that $e \not \in H_v$ but $\sigma \in H_v$ for every extremal ray $\sigma$ of $T_v(P)$ other than $e$. Similarly, choose a half space $H_w$ such that $-e \not \in H_w$ but $\sigma \in H_w$ for every extremal ray $\sigma$ of $T_w(P)$ other than $-e$. Define $Q_1=(v+T_{v}(P)) \cap (w+H_v)$, $Q_2=(v+H_w) \cap (w+T_w(P))$ and $Q_3=(v+H_w) \cap (w+H_v) \cap C$. Then $Q_1$, $Q_2$ and $Q_3$ are clean and we have $T_v(P)=T_v(Q_1)$, $T_w(Q_1)=T_w(Q_3)$, $T_v(Q_3)=T_v(Q_2)$ and $T_w(Q_2)=T_w(P)$. Using the clean case three times, we deduce the result.
\end{proof}

\begin{lemma} \label{SaturatedCones}
We have $S_{i_1, \ldots, i_d} = T_{e_{i_1}+\ldots+ e_{i_d}}(\Poly_M) \cap \ZZ^n$.
\end{lemma}

\begin{proof}
$T_{e_{i_1}+\ldots +e_{i_d}}(\Poly_M) \cap \ZZ^n$ is, essentially by definition, the saturation of $S_{i_1, \ldots, i_d} $, so our goal is to show that $S_{i_1, \ldots, i_d}$ is saturated. Let $\mathcal{T}$ be the toric variety corresponding to the polytope $\Poly_M$; our goal is to show that $\mathcal{T}$ is normal at the point corresponding to the vertex $e_{i_1} + \cdots +e_{i_d} \in \Poly_M$. But $\mathcal{T}$ is projectively normal (by \cite{White3}) so it is normal at every point. \end{proof} 

\textbf{Remark:} It is also possible to give a combinatorial proof of this lemma, which is essentially equivalent to the fact that every semigroup generated by a subset of the $A_{n+1}$ root system is saturated.

\begin{proof}[Proof of proposition \ref{NonRealizable}]
We adopt the shorthand $Bi$ and $Bj$ for $B \cup \{i \}$ and $B \cup \{ j \}$. If neither $Bi$ nor $Bj$ is a basis of $M$, then the claim is that $0 \equiv 0$, which is obvious. If $Bi$ is a basis of $M$ and $Bj$ is not then none of the edges of $\Poly_M$ coming out of vertex $v_i:=\sum_{k \in Bi} e_k$ point in direction $e_j-e_i$. Thus $h_{Bi}$ is a rational function whose denominator is not divisible by $(1-x_i/x_j)$ and we deduce that $(1-x_i/x_j)$ divides $f_{Bi}$. As $f_{Bj}=0$ in this setting, we get $
f_{Bi} \equiv f_{Bj} \mod (1-x_i/x_j)$ as desired.

 We now consider the case that $Bi$ and $Bj$ are both bases. Then $v_i:=\sum_{k \in Bi} e_k$ and $v_j:=\sum_{k \in Bj} e_k$ are both vertices of the polytope $\Poly_M$ that are joined by an edge in direction $e_i-e_j$. By lemma \ref{SaturatedCones}, $h_{Bi}=h_{T_{v_i}(\Poly_M)}$ and similarly for $Bj$. So, by lemma \ref{EdgeCancellation}, $h_{Bi}+h_{Bj}$ is a rational function whose denominator is not divisible by $(1-x_i/x_j)$. Also, each of $(1-x_j/x_i) h_{Bi}$ and  $(1-x_i/x_j) h_{Bj}$ is a rational function whose denominator is not divisible by $(1-x_i/x_j)$. 

The congruences which follow take place in the ring of rational functions whose denominator are not divisible by $(1-x_i/x_j)$.  In order to compress our equations, we adopt the shorthand $\eta(i,j)$ for $(1-x_i/x_j)$. We have
\begin{multline*}
f_{Bi} - f_{Bj} = \prod_{k \in B} \prod_{\ell \not \in Bij} \eta(\ell,k) \left( \prod_{k \in B} \eta(j,k) \prod_{\ell \not \in Bij} \eta(\ell,i) \left[ \eta(j,i) h_{Bi} \right] - \right. \\
\left. \prod_{k \in B} \eta(i,k) \prod_{\ell \not \in Bij} \eta(\ell,j) \left[ \eta(i,j) h_{Bj} \right] \right).
\end{multline*}
Notice that the terms in square brackets do not have $(1-x_i/x_j)$ in their denominator. We now concentrate on the term in the large parentheses: 
\begin{multline*}
\prod_{k \in B} \eta(j,k) \prod_{\ell \not \in Bij} \eta(\ell,i)\left[ \eta(j,i) h_{Bi} \right] - \prod_{k \in B} \eta(i,k)\prod_{\ell \not \in Bij} \eta(\ell,j)\left[ \eta(i,j) h_{Bj} \right]   \equiv \\
\prod_{k \in B} \eta(j,k) \prod_{\ell \not \in Bij} \eta(\ell,i) \left( \left[ \eta(j,i) h_{Bi} \right] -  \left[ \eta(i,j) h_{Bj} \right] \vphantom{\prod_{k \in B}} \right) \mod (1-x_i/x_j). \end{multline*}
Moving once again inside the large parentheses, 
\begin{eqnarray*}
\left[ (1-x_j/x_i) h_{Bi} \right] - \left[ (1-x_i/x_j) h_{Bj} \right]  &=& \left[ (1-x_j/x_i) h_{Bi} \right] + (x_i/x_j) \left[ (1-x_j/x_i) h_{Bi} \right] \\
& \equiv & \left[ (1-x_j/x_i) h_{Bi} \right] + \left[ (1-x_j/x_i) h_{Bi} \right] \\
& = & (1-x_i/x_j) (h_{Bi}+h_{Bj}) \equiv 0 \mod (1-x_i/x_j).
\end{eqnarray*}
So we conclude that $f_{Bi} - f_{Bj} \equiv 0 \mod (1-x_i/x_j)$, as desired.
\end{proof} 

\pagebreak

\raggedright


\thebibliography{99}

\bibitem{Brion}
M. Brion, \emph{Positivity in the Grothendieck group of complex flag varieties}, Journal of Algebra \textbf{258} (2002), no. 1, 137--159.

\bibitem{Brion2} 
M. Brion, \emph{Points entiers dans les poly\`edres convexes}, Annales Scientifiques de l'\`Ecole Normale Sup\`erieure (4) \textbf{21} (1988), no. 4, 653--663.

\bibitem{Bry}
T. Brylawski, \emph{A combinatorial model for series-parallel networks}, Transactions of the AMS \textbf{154} (1971), 1--22.

\bibitem{MaxLikeDeg}
F. Catanese, S. Hosten, A. Khetan and B. Sturmfels, \emph{The Maximum Likelihood Degree},   American Journal of Mathematics, to appear, \textbf{128} (2006).

\bibitem{Crapo}
H. Crapo, \emph{A higher invariant for matroids}, Journal of Combinatorial Theory \textbf{2} (1967),  406--417.

\bibitem{CE}
W. Cunningham and J. Edwards, \emph{A Combinatorial Decomposition Theory}, Canadian Journal of  Mathematics \textbf{32} (1980), 734--765.

\bibitem{DW}
A. Dress and W. Wenzel, \emph{Valuated Matroids}, Advances in Mathematics \textbf{93} (1992), no. 2, 214--250


\bibitem{Fult}
W. Fulton, \emph{Introduction to Toric Varieties}, Annals of Mathematics Studies Vol. 131, Princeton University Press, Princeton, 1993.

\bibitem{GGMS}
I. Gelfand, I. Goresky, R. MacPherson and V. Serganova, \emph{Combinatorial
geometries, convex polyhedra, and Schubert cells}, Advances in  Mathematics \textbf{63} (1987), 301--316.

\bibitem{GelfMac}
I. Gelfand and R. MacPherson, \emph{Geometry in Grassmannians and a generalization of the dilogarithm}, Advances in Mathematics \textbf{44} (1982), no. 3, 279 -- 312.


\bibitem{HKT} P. Hacking, S. Keel and E. Tevelev, \emph{Compactification of
  the moduli space of hyperplane arrangements}, preprint \texttt{arXiv:math.AG/0501227} 2005.

\bibitem{Ishida}
M. Ishida, \emph{Polyhedral Laurent series and Brion's equalities}, International Journal of Mathematics\textbf{1} (1990) no. 3, 251--265.

\bibitem{CQoG1} 
M. Kapranov, \emph{Chow quotients of Grassmannians I}, Advances in Soviet Mathematics \textbf{16} (1993), part 2, 29--110.

\bibitem{Kawa}
Y. Kawamata, \emph{A generalization of Kodaira-Ramanujam's vanishing theorem}, Mathematische Annalen \textbf{261} (1982), no. 1, 43--46.


\bibitem{Kleim}
S. Kleiman, \emph{The transversality of a generic translate}, Composito Mathematicae \textbf{28} (1974), 287--297.

\bibitem{KK}
B. Kostant and S. Kumar, \emph{$T$-equivariant $K$-theory of generalized flag varieties}, Journal of Differential Geometry \textbf{32} (1990), no. 2, 549--603.

\bibitem{Ezra}
E. Miller, \emph{Cohen-Macaulay quotients via irreducible resolutions}, Mathematical Research Letters \textbf{9} (2002), no. 1, 117--128.


\bibitem{MS}
E. Miller and D. Speyer, \emph{A Kleiman-Bertini theorem for sheaf tensor products}, preprint \texttt{arXiv:math.AG/0601202} 2006

\bibitem{MSturm}
E. Miller and B. Sturmfels, \emph{Combinatorial Commutative Algebra}, Graduate Studies in Mathematics Vol. 227, Springer-Verlag, New York, 2004.

\bibitem{OT} P. Orlik and H. Terao, \emph{The number of critical points of
  a product of powers of linear functions}, Inventiones Mathematicae
  \textbf{120} (1995), no. 1, 1--14.

\bibitem{Ox}
J. Oxley, \emph{On Crapo's $\beta$ invariant for matroids}, Studies in Applied Mathematics \textbf{66} (1982), no. 3, 267--277.

\bibitem{KR}
I. Rosu, \emph{Equivariant $K$-theory and equivariant cohomology}, with an appendix by A. Knutson and I. Rosu, Mathematische Zeitschrift \textbf{243} (2003), no. 3, 423--448.

\bibitem{Smirnov}
A. Smirnov, \emph{Torus schemes over a discrete valuation ring}, St. Petersburg Mathematical Journal \textbf{8} (1997), no. 4, 651--659.

\bibitem{Stan}
R. Stanley, \emph{Enumerative Combinatorics, Volume I}, Cambridge Studies in Advanced Mathematics Vol. 49, Cambridge University Press, Cambridge, 1997

\bibitem{TLS} 
D. Speyer, \emph{Tropical Linear Spaces}, preprint \texttt{ArXiv:math.CO/0410455} 2004

\bibitem{TTG}
D. Speyer and B. Sturmfels, \emph{The Tropical Grassmannian}, Advances in Geometry \textbf{4} (2004), no. 3, 389--411.

\bibitem{GrobTor}
B. Sturmfels, \emph{Gr\"obner bases of toric varieties}, Tohoku Mathematical Journal (2) \textbf{43} (1991), no. 2, 249--261.


\bibitem{Var} A. Varchenko, \emph{Critical points of the product of powers
  of linear functions and families of bases of singular vectors},
  Compositio Mathematica \textbf{97} (1995), no. 3, 385--401.
  
\bibitem{VV}
G. Vezzosi and A. Vistoli, \emph{Higher algebraic $K$-theory for actions of diagonalizable groups}, Inventiones Mathematicae \textbf{153} (2003), no. 1, 1--44.
  
\bibitem{Vieh} E. Viehweg, \emph{Vanishing theorems}, Journal f\"ur die Reine und Angewandte Mathematik \textbf{335} (1982), 1--8.  
  
\bibitem{White1}
N. White, editor,  \emph{Theory of matroids},  Encyclopedia of Mathematics and its Applications Vol. 26, Cambridge University Press, Cambridge, 1986.
  
\bibitem{White2}
N. White, editor, \emph{Matroid applications}, Encyclopedia of Mathematics and its Applications Vol. 40,  Cambridge University Press, Cambridge, 1992.  
  
 \bibitem{White3} N. White, \emph{The basis monomial ring of a matroid}, Advances in Mathematics \textbf{24} (1977), no. 3, 292--297.

\bibitem{Zieg} G. Ziegler, \emph{Lectures on Polytopes}, Graduate Texts in Mathematics Vol. 152, Springer-Verlag, New York, 1991.
\end{document}